\documentclass[12pt]{article}

\usepackage{amssymb}
\usepackage{luis}

\def\stackunder#1#2{\mathrel{\mathop{#2}\limits_{#1}}}
\def\QTR#1#2{{\csname#1\endcsname #2}}
\def\binom#1#2{{#1 \choose #2}}

\newtheorem{theorem}{Theorem}[section]

\newtheorem{corollary}[theorem]{Corollary}
\newtheorem{lemma}[theorem]{Lemma}
\newtheorem{proposition}[theorem]{Proposition}
\newtheorem{definition}[theorem]{Definition}
\newtheorem{example}[theorem]{Example}

\newtheorem{remark}[theorem]{Remark}

\begin{document}
{\LARGE
\noindent {\bf Universit\"at Bielefeld}\\[0.3cm]
\noindent Forschungszentrum\\[0.3cm]
\noindent Bielefeld-Bochum-Stochastic
}
\vspace*{5.0cm}
\begin{center}
{\Large {\bf Generalized Appell Systems}}\\[0.3cm]

Yuri G.~Kondratiev, Jos\'e L.~Silva, and Ludwig Streit\\[0.3cm]

\hspace{7cm} Nr.~729/5/96
\vspace*{6.5cm}
\begin{center}
{\bf
BiBoS\\
Universit\"at Bielefeld Postfach100131 33501 Bielefeld\\
Tel: +49 521 1065305, Fax: +49 521 1062961\\
e-mail: BiBoS @ Physik.UNI-Bielefeld.DE\\
http://www.physik.uni-bielefeld.de/bibos/start.html
}
\end{center}
\end{center}
\thispagestyle{empty}
\newpage
\title{Generalized Appell Systems}
\author{\textbf{Yuri G. Kondratiev} \\
BiBoS, Universit\"at Bielefeld, D 33615 Bielefeld, Germany \\
Institute of Mathematics, Kiev, Ukraine 
\and 
\textbf{Jos\'e Luis Silva} \\
CCM, Universidade da Madeira, P 9000 Funchal, Portugal\\
(luis@uma.pt)
\and 
\textbf{Ludwig Streit} \\
BiBoS, Universit\"at Bielefeld, D 33615 Bielefeld, Germany \\
CCM, Universidade da Madeira, P 9000 Funchal, Portugal}
\date{}
\maketitle

\begin{abstract}
We give a general approach to infinite dimensional non-Gaussian analysis
which generalizes the work \cite{KSWY95}. For given measure we construct a
family of biorthogonal systems. We study their properties and their Gel'fand
triples that they generate. As an example we consider the measures of
Poisson type.
\end{abstract}
\renewcommand{\thefootnote}{}
\footnote{
\noindent Published in \textit{Methods of Functional Analysis and Topology}, 
Vol.~\textbf{3} no.~3, 1997.

Available on WWW: http://www.physik.uni-bielefeld.de/bibos/start.html}
\thispagestyle{empty}

\newpage 
\setcounter{page}{1}
\renewcommand{\thefootnote}{\arabic{footnote}}

\tableofcontents

\section{Introduction}

Non-Gaussian analysis was already introduced in \cite{AKS93} for smooth
proba\-bility measure on infinite dimensional linear spaces, using
biorthogonal decomposition as a natural extension of the chaos decomposition
that is well known in Gaussian analysis. This biorthogonal ``Appell'' system
has been constructed for smooth measures by Yu. L. Daletskii \cite{Da91}.
For a detailed description of its use in infinite dimensional analysis and
for the proof of the results which were announced in \cite{AKS93} we refer
to \cite{ADKS96} which was based on quasi-invariance of the measures and
smoothness of the logarithmic derivatives.

Kondratiev et al. \cite{KSWY95} considered the case of non-degenerate
measures on the dual of a nuclear space with analytic characteristic
functionals and no further conditions such as quasi-invariance of the
measure or smoothness of the logarithmic derivative was required. In this
case the important example of Poisson noise is now accessible. Again for a
given measure $\mu $ with analytic Laplace transform \cite{KSWY95} construct
an Appell biorthogonal system \textrm{\textbf{A}}$^{\mu }$ as a pair (%
\textrm{\textbf{P}}$^{\mu }{,}$\textrm{\textbf{Q}}$^{\mu }$) of Appell
polynomials \textrm{\textbf{P}}$^{\mu }$ and a canonical system of
generalized functions \textrm{\textbf{Q}}$^{\mu }$, properly associated to
the measure $\mu .$ Hence within this framework they obtained:

\begin{itemize}
\item  explicit description of the test function space introduced in \cite
{ADKS96};

\item  the test functions space is identical for all measures that they
consider;

\item  characterization theorems for generalized as well as test functions
was obtained analogously as in Gaussian analysis, see \cite{KLPSW96} for
more references;

\item  extension of the Wick product and the corresponding Wick calculus 
\cite{KLS96} as well as full description of positive distributions (as
measures).\smallskip 
\end{itemize}

\textbf{Aim of the present work}. As in \cite{KSWY95} we consider the case
of non-degenerate measures on the dual of a nuclear space with analytic
Laplace transform but instead of the $\mu $\textbf{-exponential} $e_{\mu
}(\cdot ,\cdot )$ we use the \textbf{generalized}\textit{\ }$\mu $\textbf{%
-exponential} $e_{\mu }^{\alpha }(\cdot ,\cdot )$ where $\alpha $ is a
holomorphic function $\alpha $ on $\mathcal{N}_{\QTR{mathbb}{C}}$ which is
invertible in a neighborhood of zero, i.e., $\alpha \in \mathrm{Hol}_{0}(%
\mathcal{N}_{\QTR{mathbb}{C}},\mathcal{N}_{\QTR{mathbb}{C}}).$ Hence using $%
e_{\mu }^{\alpha }(\cdot ,\cdot )$ we construct an generalized Appell
orthogonal system \textrm{\textbf{A}}$^{\mu ,\alpha }$ as a pair (\textrm{%
\textbf{P}}$^{\mu ,\alpha }{,}$\textrm{\textbf{Q}}$^{\mu ,\alpha }$) of
generalized Appell polynomials \textrm{\textbf{P}}$^{\mu ,\alpha }$ and a
system of genera\-lized functions \textrm{\textbf{Q}}$^{\mu ,\alpha }.$

\noindent \textbf{Central results}. In the above framework

\begin{itemize}
\item  we obtain an explicit description of the test function space
introduced in \cite{ADKS96};

\item  the spaces of test functions turns out to be the same for all $\alpha
\in \mathrm{Hol}_{0}(\mathcal{N}_{\QTR{mathbb}{C}},\mathcal{N}_{\QTR{mathbb}{%
C}})$ and for all measures that we consider;

\item  characterization theorems for generalized as well as test functions
are obtained analogously as in the Gaussian case;

\item  the spaces of distributions for a fixed measure $\mu $ are again
identical for all function $\alpha $ in the above conditions;

\item  the well known Wick product and the corresponding Wick calculus \cite
{KLS96} extends rather directly;

\item  in the important case of Poisson white noise a special choice of $%
\alpha $ produces the orthogonal system of Charlier polynomials, see Example 
\ref{1eq32}.
\end{itemize}

\section{General theory\label{1eq8}}

\subsection{Some facts on nuclear triples}

We start with a real separable Hilbert space $\mathcal{H}$ with inner
product $(\cdot ,\cdot )$ and norm $\left| \cdot \right| $. For a given
separable nuclear space $\mathcal{N}$ densely topologically embedded in $%
\mathcal{H}$ we can construct the nuclear triple
\[
\mathcal{N}\subset \mathcal{H}\subset \mathcal{N}^{\prime }.
\]
The dual pairing $\langle \cdot ,\cdot \rangle $ of $\mathcal{N}^{\prime }$
and $\mathcal{N}$ then is realized as an extension of the inner product in $%
\mathcal{H}$%
\[
\langle f,\,\xi \rangle =(f,\xi )\quad f\in \mathcal{H},\ \xi \in \mathcal{N}%
.
\]
Instead of reproducing the abstract definition of nuclear spaces (see e.g., 
\cite{S71}) we give a complete (and convenient) characterization in terms of
projective limits of decreasing chains of Hilbert spaces $\mathcal{H}_{p}$, $%
p\in \QTR{mathbb}{N}$.

\begin{theorem}
The nuclear Fr\'{e}chet space $\mathcal{N}$ can be represented as 
\[
\mathcal{N}=\bigcap_{p\in \QTR{mathbb}{N}}\mathcal{H}_{p},
\]
where $\{\mathcal{H}_{p},\;p\in \QTR{mathbb}{N}{\}}$ is a family of Hilbert
spaces such that for all $p_{1},p_{2}\in \QTR{mathbb}{N}$ there exists $p\in 
\QTR{mathbb}{N}$ such that the embeddings $\mathcal{H}_{p}\hookrightarrow 
\mathcal{H}_{p_{1}}$, $\mathcal{H}_{p}\hookrightarrow \mathcal{H}_{p_{2}}$
are of Hilbert-Schmidt class. The topology of $\mathcal{N}$ is given by the
projective limit topology, i.e., the coarsest topology on $\mathcal{N}$ such
that the canonical embeddings $\mathcal{N}\hookrightarrow \mathcal{H}_{p}$
are continuous for all $p\in \QTR{mathbb}{N}$.
\end{theorem}

The Hilbert norms on $\mathcal{H}_{p}$ are denoted by $\left| \cdot \right|
_{p}$. Without loss of generality we always suppose that $\forall p\in 
\QTR{mathbb}{N}$, $\forall \xi \in \mathcal{N}:$ $\left| \xi \right| \leq
\left| \xi \right| _{p}$ and that the system of norms is ordered, i.e., $%
\left| \cdot \right| _{p}$ $\leq \left| \cdot \right| _{q}$ if $p<q$. By
general duality theory the dual space $\mathcal{N}^{\prime }\,$ can be
written as 
\[
\mathcal{N}^{\prime }=\bigcup\limits_{p\in \QTR{mathbb}{N}}\mathcal{H}_{-p}.
\]
with inductive limit topology $\tau _{ind}$ by using the dual family of
spaces $\{\mathcal{H}_{-p}:=\mathcal{H}_{p}^{\prime },\ p\in \QTR{mathbb}{N}{%
\}}$. The inductive limit topology (w.r.t. this family) is the finest
topology on $\mathcal{N}^{\prime }$ such that the embeddings $\mathcal{H}%
_{-p}\hookrightarrow \mathcal{N}^{\prime }$ are continuous for all $p\in 
\QTR{mathbb}{N}$. It is convenient to denote the norm on $\mathcal{H}_{-p}$
by $\left| \cdot \right| _{-p}$. Let us mention that in our setting the
topology $\tau _{ind}$ coincides with the Mackey topology $\tau (\mathcal{N}%
^{\prime },\mathcal{N})$ and the strong topology $\beta (\mathcal{N}^{\prime
},\mathcal{N})$, see e.g., \cite[Appendix 5]{HKPS93}.

Further we want to introduce the notion of tensor power of a nuclear space.
The simplest way to do this is to start from usual tensor powers $\mathcal{H}%
_{p}^{\otimes n}$,$\ n\in \QTR{mathbb}{N}$ of Hilbert spaces. Since there is
no danger of confusion we will preserve the notation $\left| \cdot \right|
_{p}$ and $\left| \cdot \right| _{-p}$ for the norms on $\mathcal{H}%
_{p}^{\otimes n}$ and $\mathcal{H}_{-p}^{\otimes n}$ respe\-ctively. Using
the definition 
\[
\mathcal{N}^{\otimes n}:=\ \stackunder{p\in \QTR{mathbb}{N}}{\mathrm{%
pr\,\lim }}\mathcal{H}_{p}^{\otimes n},
\]
one can prove \cite{S71} that $\mathcal{N}^{\otimes n}$ is a nuclear space
which is called the n-th tensor power of $\mathcal{N}$.

The dual space of $\mathcal{N}^{\otimes n}$ can be written 
\[
\mathcal{N}^{\prime \otimes n}=\ \stackunder{p\in \QTR{mathbb}{N}}{\mathrm{%
ind\,}\lim }\mathcal{H}_{-p}^{\otimes n}.
\]

We also want to introduce the (Boson or symmetric) Fock space $\Gamma (%
\mathcal{H})$ of $\mathcal{H}$ by 
\[
\Gamma (\mathcal{H})=\bigoplus_{n=0}^{\infty }\mathcal{H}_{\QTR{mathbb}{C}}^{%
\widehat{\otimes }n}
\]
with the convention $\mathcal{H}_{\QTR{mathbb}{C}}^{\widehat{\otimes }0}:=%
\QTR{mathbb}{C}$ and the Hilbert norm 
\[
\left\| \vec{\varphi}\right\| _{\Gamma (\mathcal{H})}^{2}=\sum_{n=0}^{\infty
}n!\left| \varphi ^{(n)}\right| ^{2},\;\vec{\varphi}=\left\{ \varphi
^{(n)}\,|\,n\in \QTR{mathbb}{N}_{0}\right\} \in \Gamma (\mathcal{H}).
\]

\subsection{Holomorphy on locally convex spaces \label{1eq9}}

We shall collect some facts from the theory of holomorphic functions in
locally convex topological vector spaces $\mathcal{E}$ (over the complex
field $\QTR{mathbb}{C}$), see e.g., \cite{Di81}. Let $\mathcal{L}(\mathcal{E}%
^{n})$ be the space of n-linear mappings from $\mathcal{E}^{n}$ into $%
\QTR{mathbb}{C}$ and $\mathcal{L}_{s}(\mathcal{E}^{n})$ the subspace of
symmetric n-linear forms. Also let $P^{n}(\mathcal{E)}$ denote the
n-homogeneous polynomials on $\mathcal{E}$. There is a linear bijection $%
\mathcal{L}_{s}(\mathcal{E}^{n})\ni A\longleftrightarrow \widehat{A}\in
P^{n}(\mathcal{E)}$. Now let $\mathcal{U}\subset \mathcal{E}$ be open and
consider a function $G:\mathcal{U}\rightarrow \QTR{mathbb}{C}$. $G$ is said
to be \textbf{G-holomorphic} if for all $\theta _{0}\in \mathcal{U}$ and for
all $\theta \in \mathcal{E}$ the mapping from $\QTR{mathbb}{C}$ to $%
\QTR{mathbb}{C}{:}$ $\lambda \mapsto G(\theta _{0}+\lambda \theta )$ is
holomorphic in some neighborhood of zero in $\QTR{mathbb}{C}$. If $G$ is
G-holomorphic then there exists for every $\eta \in \mathcal{U}$ a sequence
of homogeneous polynomials $\frac{1}{n!}\widehat{d^{n}G(\eta )}$ such that 
\[
G(\theta +\eta )=\sum\limits_{n=0}^{\infty }\frac{1}{n!}\widehat{d^{n}G(\eta
)}(\theta )
\]
for all $\theta $ from some open neighborhood $\mathcal{V}$ of zero. $G$ is
said to be \textbf{holomorphic}, if for all $\eta $ in $\mathcal{U}$ there
exists an open neighborhood $\mathcal{V}$ of zero such that 
\[
\sum_{n=0}^{\infty }\frac{1}{n!}\widehat{d^{n}G(\eta )}(\theta )
\]
converges uniformly on $\mathcal{V}$ to a continuous function. Of course, $%
\widehat{d^{n}G(\eta )}(\theta )$ is the n-th partial derivative of $G$ at $%
\eta $ in direction $\theta $. We say that $G$ is holomorphic at $\theta _{0}
$ if there is an open set $\mathcal{U}$ containing $\theta _{0}$ such that $G
$ is holomorphic on $\mathcal{U}$. The following Proposition can be found
e.g., in \cite{Di81}.

\begin{proposition}
\label{1eq10}$G$ is holomorphic if and only if it is G-holomorphic and
locally bounded.
\end{proposition}

\noindent Let us explicitly consider a function holomorphic at the point $%
0\in \mathcal{E}=\mathcal{N}_{\QTR{mathbb}{C}}$, then

1) there exist $p$ and $\varepsilon >0$ such that for all $\xi _{0}\in 
\mathcal{N}_{\QTR{mathbb}{C}}$ with $\left| \xi _{0}\right| _{p}\leq
\varepsilon $ and for all $\xi \in \mathcal{N}_{\QTR{mathbb}{C}}$ the
function of one complex variable $\lambda \mapsto G(\xi _{0}+\lambda \xi )$
is holomorphic at $0\in \QTR{mathbb}{C}$, and

2) there exists $c>0$ such that for all $\xi \in \mathcal{N}_{\QTR{mathbb}{C}%
}$ with $\left| \xi \right| _{p}\leq \varepsilon $ : $\left| G(\xi )\right|
\leq c$.

\noindent As we do not want to discern between different restrictions of one
function, we consider germs of holomorphic functions, i.e., we identify $F$
and $G$ if there exists an open neighborhood $\mathcal{U}:0\in \mathcal{U}%
\subset \mathcal{N}_{\QTR{mathbb}{C}}$ such that $F(\xi )=G(\xi )$ for all $%
\xi \in \mathcal{U}$. Thus we define $\mathrm{Hol}_{0}(\mathcal{N}_{%
\QTR{mathbb}{C}})$ as the algebra of germs of functions holomorphic at zero
equipped with the inductive topology given by the following family of norms 
\[
\mathrm{n}_{p,l,\infty }(G)=\sup_{\left| \theta \right| _{p}\leq
2^{-l}}\left| G(\theta )\right| ,\quad p,l\in \QTR{mathbb}{N}.
\]

For later use we need the space $\mathrm{Hol}_{0}(\mathcal{N}_{\QTR{mathbb}{C%
}},\mathcal{N}_{\QTR{mathbb}{C}})$ of holomorphic functions from $\mathcal{N}%
_{\QTR{mathbb}{C}}$ to $\mathcal{N}_{\QTR{mathbb}{C}}$. Let $\mathcal{U}%
\subset \mathcal{N}_{\QTR{mathbb}{C}}$ be open and consider a function $%
\alpha :\mathcal{U}\rightarrow \mathcal{N}_{\QTR{mathbb}{C}}$. $\alpha $ is
said to be holomorphic at $0\in \mathcal{N}_{\QTR{mathbb}{C}}$ iff

\begin{enumerate}
\item  it is G-holomorphic; i.e., there exist $p$ and $\epsilon >0$ such
that for all $\xi _{0}\in \mathcal{N}_{\QTR{mathbb}{C}}$ with $\left| \xi
_{0}\right| _{p}\leq \epsilon $ and for all $\xi \in \mathcal{N}_{%
\QTR{mathbb}{C}}$ the function of one complex variable $\lambda \mapsto
\alpha (\xi _{0}+\lambda \xi )$ is holomorphic at $0\in \QTR{mathbb}{C}${;}

\item  $\alpha $ is locally bounded, i.e., for all $p\in \QTR{mathbb}{N}$
there exist $C_{p}>0$ such that $\forall \eta \in A$ with $\left| \eta
\right| _{p}\leq C_{p}$ then $\forall p^{\prime }\in \QTR{mathbb}{N}$ there
exist $C_{p^{\prime }}$ such that $\forall \eta \in A$ $\left| \alpha (\eta
)\right| _{p^{\prime }}\leq C_{p^{\prime }}$, where $A$ is an bounded set in 
$\mathcal{N}_{\QTR{mathbb}{C}}$.
\end{enumerate}

\noindent If $\alpha $ is holomorphic at $0\in \mathcal{N}_{\QTR{mathbb}{C}}$%
, then for every $\eta \in \mathcal{U}$ there exists a sequence of
homogeneous polynomials $\frac{1}{n!}\widehat{d^{n}\alpha (\eta )}$ such
that 
\[
\theta \longmapsto \sum_{n=0}^{\infty }\frac{1}{n!}\widehat{d^{n}\alpha
(\eta )}\left( \theta \right) 
\]
converges and define a continuous function on some neighborhood of zero.

Let use now introduce spaces of entire functions which will be useful later.
Let $\mathcal{E}_{2^{-l}}^{k}(\mathcal{H}_{-p,\QTR{mathbb}{C}})$ denote the
set of all entire functions on $\mathcal{H}_{-p,\QTR{mathbb}{C}}$ of growth $%
k\in [1,2]$ and type $2^{-l},\ p,l\in \QTR{mathbb}{Z}$. This is a linear
space with norm 
\[
\mathrm{n}_{p,l,k}(\varphi )=\sup_{z\in \mathcal{H}_{-p,\QTR{mathbb}{C}%
}}\left| \varphi (z)\right| \exp \left( -2^{-l}|z|_{-p}^{k}\right)
,\;\varphi \in \mathcal{E}_{2^{-l}}^{k}(\mathcal{H}_{-p,\QTR{mathbb}{C}}).
\]
The space of entire functions on $\mathcal{N}_{\QTR{mathbb}{C}}^{\prime }$
of growth $k$ and minimal type is naturally introduced by 
\[
\mathcal{E}_{\min }^{k}(\mathcal{N}_{\QTR{mathbb}{C}}^{\prime }):=\ 
\stackunder{p,l\in \QTR{mathbb}{N}}{\mathrm{pr\,\lim \,}}\mathcal{E}%
_{2^{-l}}^{k}(\mathcal{H}_{-p,\QTR{mathbb}{C}}),\ 
\]
see e.g., \cite{Ko80a}, \cite{BK88}, \cite{Ou91}. We will also need the
space of entire functions on $\mathcal{N}_{\QTR{mathbb}{C}}$ of growth $k$
and finite type: 
\[
\mathcal{E}_{\max }^{k}(\mathcal{N}_{\QTR{mathbb}{C}}):=\ \stackunder{p,l\in 
\QTR{mathbb}{N}}{\mathrm{ind\,\lim }}\mathcal{E}_{2^{l}}^{k}(\mathcal{H}_{p,%
\QTR{mathbb}{C}}).
\]

\subsection{Measures on linear topological spaces}

To introduce probability measures on the vector space $\mathcal{N}^{\prime }$%
, we consider $\mathcal{C}_{\sigma }(\mathcal{N}^{\prime })$ the $\sigma $%
-algebra generated by cylinder sets on $\mathcal{N}^{\prime }$, which
coincides with the Borel $\sigma $-algebras $\mathcal{B}_{\sigma }(\mathcal{N%
}^{\prime })$ and $\mathcal{B}_{\beta }(\mathcal{N}^{\prime })$ generated by
the weak and strong topology on $\mathcal{N}^{\prime }$, respectively. Thus
we will consider this $\sigma $-algebra as the \textbf{natural} $\sigma $%
-algebra on $\mathcal{N}^{\prime }$. Detailed definitions of the above
notions and proofs of the mentioned relations can be found in e.g., \cite
{BK88}.

We will restrict our investigations to a special class of measures $\mu $ on 
$\mathcal{C}_{\sigma }(\mathcal{N}^{\prime })$ which satisfy two additional
assumptions. The first one concerns some analyticity of the Laplace
transformation 
\[
l_{\mu }\left( \theta \right) :=L_{\mu }1\left( \theta \right) =\int_{%
\mathcal{N}^{\prime }}\exp \left\langle x,\theta \right\rangle d\mu (x)=:%
\QTR{mathbb}{E}_{\mu }\left( \exp \left\langle \cdot ,\theta \right\rangle
\right) ,\quad \theta \in \mathcal{N}_{\QTR{mathbb}{C}}\,.
\]
Here we also have introduced the convenient notion of expectation 
$\QTR{mathbb}{E}_{\mu }$ of a $\mu $-integrable function. \bigskip 

\noindent \textbf{Assumption 1 }The measure $\mu $ has an analytic Laplace
transform in a neighborhood of zero. That means there exists an open
neighborhood $\mathcal{U}\subset \mathcal{N}_{\QTR{mathbb}{C}}$ of zero,
such that $l_{\mu }$ is holomorphic on $\mathcal{U}$, i.e., $l_{\mu }\in 
\mathrm{Hol}_{0}(\mathcal{N}_{\QTR{mathbb}{C}})$. This class of \textbf{%
analytic measures} is denoted by $\mathcal{M}_{a}(\mathcal{N}^{\prime })$%
.\bigskip 

\noindent An equivalent description of analytic measures is given by the
following lemma and the proof can be founded in \cite{KSW95}.

\begin{lemma}
\label{1eq11}The following statements\thinspace are equivalent
\end{lemma}

\textbf{1)}\quad $\mu \in \mathcal{M}_{a}(\mathcal{N}^{\prime })$;\smallskip 

\textbf{2)}$\quad \displaystyle\exists p_{\mu }\in \QTR{mathbb}{N},\quad
\exists C>0:\;\left| \int_{\mathcal{N}^{\prime }}\langle x,\theta \rangle
^{n}d\mu (x)\right| \leq n!\,C^{n}\left| \theta \right| _{p_{\mu
}}^{n},\;\theta \in \mathcal{H}_{p_{\mu },\QTR{mathbb}{C}};$\smallskip 

\textbf{3)}$\quad \displaystyle\exists p_{\mu }^{\prime }\in \QTR{mathbb}{N}%
,\quad \exists \varepsilon _{\mu }>0:\;\int_{\mathcal{N}^{\prime }}\exp
(\varepsilon _{\mu }\left| x\right| _{-p_{\mu }^{\prime }})d\mu (x)<\infty .$%
\smallskip 

\bigskip \ For $\mu \in \mathcal{M}_{a}(\mathcal{N}^{\prime })$ the estimate
in statement 2 of the above lemma allows to define the moment kernels $%
M_{n}^{\mu }\in \mathcal{N}^{\prime \widehat{\otimes }n}$. This is done by
extending the above estimate by a simple polarization argument and applying
the kernel theorem. The kernels are determined by 
\begin{equation}
l_{\mu }(\theta )=\sum_{n=0}^{\infty }\frac{1}{n!}\left\langle M_{n}^{\mu
},\theta ^{\otimes n}\right\rangle   \label{1eq12}
\end{equation}
or equivalently 
\[
\left\langle M_{n}^{\mu },\theta _{1}\widehat{\otimes }\ldots \widehat{%
\otimes }\theta _{n}\right\rangle =\left. \frac{\partial ^{n}}{\partial
t_{1}\ldots \partial t_{n}}l_{\mu }\left( t_{1}\theta _{1}+\ldots
+t_{n}\theta _{n}\right) \right| _{t_{1}=\ldots =t_{n}=0}.
\]
Moreover, if $p>p_{\mu }$ is such that the embedding $i_{p,p_{\mu }}:%
\mathcal{H}_{p}\hookrightarrow \mathcal{H}_{p_{\mu }}$ is Hilbert-Schmidt
then 
\begin{equation}
\left| M_{n}^{\mu }\right| _{-p}\leq \left( nC\left\| i_{p,p_{\mu }}\right\|
_{HS}\right) ^{n}\leq n!\left( eC\left\| i_{p,p_{\mu }}\right\| _{HS}\right)
^{n}.  \label{1eq13}
\end{equation}

\begin{definition}
\label{1eq14}A function $\varphi :\mathcal{N}^{\prime }\rightarrow 
\QTR{mathbb}{C}$ of the form 
\[
\varphi (x)=\sum_{n=0}^{N}\langle x^{\otimes n},\varphi ^{(n)}\rangle
,\;x\in \mathcal{N}^{\prime },\;N\in \QTR{mathbb}{N}{,}
\]
is called a \emph{\textbf{continuous polynomial }}$($short $\varphi \in 
\mathcal{P}(\mathcal{N}^{\prime }))$ iff $\varphi ^{(n)}\in \mathcal{N}_{%
\QTR{mathbb}{C}}^{\widehat{\otimes }n}$, $\forall n\in \QTR{mathbb}{N}_{0}:=%
\QTR{mathbb}{N}\cup \{0\}$.
\end{definition}

Now we are ready to formulate the second assumption on $\mu $:\bigskip

\noindent \textbf{Assumption 2 }For all $\varphi \in \mathcal{P}(\mathcal{N}%
^{\prime })$ with $\varphi =0$ $\mu $-almost everywhere we have $\varphi
\equiv 0$. In the following a measure with this property will be called 
\textbf{non-degenerate}.\bigskip

\noindent \textbf{Note:} Assumption 2 is equivalent to:

\noindent Let $\varphi \in \mathcal{P}(\mathcal{N}^{\prime })$ with $%
\int_{A}\varphi d\mu =0$ for all $A\in \mathcal{C}_{\sigma }(\mathcal{N}%
^{\prime })$ then $\varphi \equiv 0$.

\noindent A sufficient condition can be obtained by regarding admissible
shifts of the measure $\mu $. If $\mu (\cdot +\xi )$ is absolutely
continuous with respect to $\mu $ for all $\xi \in \mathcal{N}$, i.e., there
exists the Radon-Nikodym derivative 
\[
\rho _{\mu }\left( \xi ,x\right) =\frac{d\mu \left( x+\xi \right) }{d\mu
\left( x\right) }\in L^{1}\left( \mathcal{N}^{\prime },\mu \right) ,\;x\in 
\mathcal{N}^{\prime },
\]
then we say that $\mu $ is $\mathcal{N}$\textbf{-quasi-invariant}\textit{\ }%
see e.g., \cite{GV68}, \cite{Sk74}. This is sufficient to ensure Assumption
2, see e.g., \cite{KoTs91}, \cite{BK88}.

\section{The Appell system\label{1eq43}}

The space $\mathcal{P}(\mathcal{N}^{\prime })$ may be equipped with various
different topologies, but there exists a natural one such that $\mathcal{P}(%
\mathcal{N}^{\prime })$ becomes isomorphic to the topological direct sum of
tensor powers $\mathcal{N}_{\QTR{mathbb}{C}}^{\widehat{\otimes }n}$ see
e.g., \cite[Chap. II 6.1, Chap. II 7.4]{S71} 
\[
\mathcal{P}(\mathcal{N}^{\prime })\simeq \bigoplus_{n=0}^{\infty }\mathcal{N}%
_{\QTR{mathbb}{C}}^{\widehat{\otimes }n} 
\]
via 
\[
\varphi \left( x\right) =\sum_{n=0}^{\infty }\left\langle x^{\otimes
n},\varphi ^{\left( n\right) }\right\rangle \longleftrightarrow \vec{\varphi}%
=\left\{ \varphi ^{\left( n\right) }\,|\,n\in \QTR{mathbb}{N}_{0}\right\} . 
\]
Note that only a finite number of $\varphi ^{\left( n\right) }$ is a
non-zero. The notion of convergence of sequences in this topology on $%
\mathcal{P}(\mathcal{N}^{\prime })$ is the following: for $\varphi \in 
\mathcal{P}(\mathcal{N}^{\prime })$, such that 
\[
\varphi \left( x\right) =\sum_{n=0}^{N(\varphi )}\left\langle x^{\widehat{%
\otimes }n},\varphi ^{\left( n\right) }\right\rangle 
\]
let $p_{n}:\mathcal{P}(\mathcal{N}^{\prime })\rightarrow \mathcal{N}_{%
\QTR{mathbb}{C}}^{\widehat{\otimes }n}$ denote the mapping $p_{n}$ defined
by $p_{n}(\varphi ):=\varphi ^{\left( n\right) }.$ A sequence $\left\{
\varphi _{j},\,j\in \QTR{mathbb}{N}\right\} $ of smooth polynomials converge
to $\varphi \in \mathcal{P}(\mathcal{N}^{\prime })$ iff the $N(\varphi _{j})$
are bounded and $p_{n}\varphi _{j}\stackunder{j\rightarrow \infty }{%
\longrightarrow }p_{n}\varphi $ in $\mathcal{N}_{\QTR{mathbb}{C}}^{\widehat{%
\otimes }n}$ for all $n\in \QTR{mathbb}{N}{.}$

Now we can introduce the dual space $\mathcal{P}_\mu ^{\prime }(\mathcal{N}%
^{\prime })$ of $\mathcal{P}(\mathcal{N}^{\prime })$ with respect to $%
L^2(\mu ).$ As a result we have constructed the triple 
\[
\mathcal{P}(\mathcal{N}^{\prime })\subset L^2(\mu )\subset \mathcal{P}_\mu
^{\prime }(\mathcal{N}^{\prime }) 
\]
The (bilinear) dual pairing $\left\langle \!\left\langle \cdot ,\cdot
\right\rangle \!\right\rangle _\mu $ between $\mathcal{P}_\mu ^{\prime }(%
\mathcal{N}^{\prime })$ and $\mathcal{P}(\mathcal{N}^{\prime })$ is
connected to the (sesquilinear) inner product on $L^2(\mu )$ by 
\[
\left\langle \!\left\langle \varphi ,\psi \right\rangle \!\right\rangle _\mu
=(\varphi ,\overline{\psi })_{L^2(\mu )}\,,\quad \varphi \in L^2\left( \mu
\right) ,\;\psi \in \mathcal{P}\left( \mathcal{N}^{\prime }\right) . 
\]

\subsection{\textrm{\textbf{P}}$^{\mu }${-}system}

Because of the holomorphy of $l_{\mu }$ and since that $l_{\mu }(0)=1$,
there exists a neighborhood of zero 
\[
\mathcal{U}_{0}=\left\{ \theta \in \mathcal{N}_{\QTR{mathbb}{C}}\ |\
2^{q_{0}}\left| \theta \right| _{p_{0}}<1\right\} 
\]
$p_{0},q_{0}\in \QTR{mathbb}{N}${,} $p_{0}\geq p_{\mu }^{\prime }$, $%
2^{-q_{0}}\leq \varepsilon _{\mu }$ ($p_{\mu }^{\prime },\varepsilon _{\mu }$
from Lemma \ref{1eq11}) such that $l_{\mu }(\theta )\neq 0$ for $%
\theta \in \mathcal{U}_{0}$ and the normalized $\mu $-exponential 
\begin{equation}
e_{\mu }\left( \theta ;z\right) :=\frac{\exp \left\langle z,\theta
\right\rangle }{l_{\mu }(\theta )}\quad \mathrm{for}\;\theta \in \mathcal{U}%
_{0},\quad z\in \mathcal{N}_{\QTR{mathbb}{C}}^{\prime }\,,  \label{1eq44}
\end{equation}
is well defined. We use the holomorphy of $\theta \mapsto e_{\mu }(\theta
;z) $ to expand it in a power series in $\theta $ similar to the case
corresponding to the construction of one dimensional Appell polynomials \cite
{Bo76}. We have in analogy to \cite{AKS93}, \cite{ADKS96} 
\[
e_{\mu }\left( \theta ;z\right) =\sum_{n=0}^{\infty }\frac{1}{n!}\widehat{%
d^{n}e_{\mu }\left( 0,z\right) }\left( \theta \right) , 
\]
where $\widehat{d^{n}e_{\mu }\left( 0,z\right) }$ is an n-homogeneous form
polynomial. But since $e_{\mu }\left( \theta ;z\right) $ is not only
G-holomorphic but holomorphic we know that $\theta \mapsto e_{\mu }\left(
\theta ;z\right) $ is also locally bounded. Thus Cauchy's inequality for
Taylor series \cite{Di81} may be applied, $\rho \leq 2^{-q_{0}}$, $p\geq
p_{0}$%
\begin{eqnarray}
\left| \frac{1}{n!}\widehat{d^{n}e_{\mu }\left( 0,z\right) }\left( \theta
\right) \right| &\leq &\frac{1}{\rho ^{n}}\sup_{\left| \theta \right|
_{p}=\rho }\left| e_{\mu }\left( \theta ;z\right) \right| \left| \theta
\right| _{p}^{n}  \nonumber \\
&\leq &\frac{1}{\rho ^{n}}\sup_{\left| \theta \right| _{p}=\rho }\frac{1}{%
l_{\mu }\left( \theta \right) }\exp \left( \rho \left| z\right| _{-p}\right)
\left| \theta \right| _{p}^{n}  \label{1eq45}
\end{eqnarray}
if $z\in \mathcal{H}_{-p,\QTR{mathbb}{C}}.$ This inequality extends by
polarization \cite{Di81} to an estimate sufficient for the kernel theorem.
Thus we have a representation 
\[
\widehat{d^{n}e_{\mu }\left( 0,z\right) }\left( \theta \right) =\left\langle
P_{n}^{\mu }(z),\theta ^{\otimes n}\right\rangle , 
\]
where $P_{n}^{\mu }(z)\in \mathcal{N}_{\QTR{mathbb}{C}}^{\prime \widehat{%
\otimes }n}.$ The kernel theorem really gives a little more: $P_{n}^{\mu
}(z)\in \mathcal{H}_{-p^{\prime },\QTR{mathbb}{C}}^{\widehat{\otimes }n}$
for any $p^{\prime }$ ($>p\geq p_{0}$) such that the embedding operator 
\[
i_{p^{\prime },p}:\mathcal{H}_{p^{\prime },\QTR{mathbb}{C}}\hookrightarrow 
\mathcal{H}_{p,\QTR{mathbb}{C}} 
\]
is Hilbert-Schmidt. Thus we have 
\begin{equation}
e_{\mu }(\theta ;z)=\sum_{n=0}^{\infty }\frac{1}{n!}\left\langle P_{n}^{\mu
}(z),\theta ^{\otimes n}\right\rangle \quad \mathrm{for}\;\theta \in 
\mathcal{U}_{0},\;z\in \mathcal{N}_{\QTR{mathbb}{C}}^{\prime }.
\label{1eq46}
\end{equation}
We will also use the notation 
\[
P_{n}^{\mu }(\varphi ^{(n)})(\cdot ):=\left\langle P_{n}^{\mu }(\cdot
),\varphi ^{(n)}\right\rangle ,\quad \varphi ^{(n)}\in \mathcal{N}_{%
\QTR{mathbb}{C}}^{\widehat{\otimes }n},\quad n\in \QTR{mathbb}{N}, 
\]
which is called \textbf{Appell polynomial}. Thus for any measure satisfying
Assumption 1 we have defined the \textrm{\textbf{P}}$^{\mu }$-system

\[
\mathrm{\mathbf{P}}^{\mu }=\left\{ \left\langle P_{n}^{\mu }(\cdot ),\varphi
^{(n)}\right\rangle \ |\ \varphi ^{(n)}\in \mathcal{N}_{\QTR{mathbb}{C}}^{%
\widehat{\otimes }n},\ n\in \QTR{mathbb}{N}_{0}\right\} . 
\]

The following proposition collects some properties of the polynomials $%
P_n^\mu \left( z\right) $, (for the proof we refer to \cite{KSWY95}).

\begin{proposition}
For $x,y\in \mathcal{N}^{\prime }$,$\ n\in \QTR{mathbb}{N}$ the following
holds\\[6mm]
(P1)\\[-11mm]
\begin{equation}
P_{n}^{\mu }(x)=\sum_{k=0}^{n}\binom{n}{k}x^{\otimes k}\widehat{\otimes }%
P_{n-k}^{\mu }(0).  \label{1eq47}
\end{equation}
\\[5mm]
(P2)\\[-11mm]
\begin{equation}
x^{\otimes n}=\sum_{k=0}^{n}\binom{n}{k}P_{k}^{\mu }(x)\widehat{\otimes }%
M_{n-k}^{\mu }.  \label{1eq48}
\end{equation}
\\[5mm]
(P3)\\[-11mm]
\[
P_{n}^{\mu }(x+y)=\sum_{k+l+m=n}\frac{n!}{k!\,l!\,m!}P_{k}^{\mu }(x)\widehat{%
\otimes }P_{l}^{\mu }(y)\widehat{\otimes }M_{m}^{\mu }
\]
\begin{equation}
=\sum_{k=0}^{n}\binom{n}{k}P_{k}^{\mu }(x)\widehat{\otimes }y^{\otimes
(n-k)}.  \label{1eq49}
\end{equation}
(P4)\ Further we observe 
\begin{equation}
\QTR{mathbb}{E}_{\mu }(\langle P_{m}^{\mu }(\cdot ),\varphi ^{(m)}\rangle
)=0\quad \mathrm{for\;}m\neq 0\ ,\varphi ^{(m)}\in \mathcal{N}_{\QTR{mathbb}{%
C}}^{\widehat{\otimes }m}.  \label{1eq50}
\end{equation}
(P5) For all $p>p_{0}$ such that the embedding $\mathcal{H}%
_{p}\hookrightarrow \mathcal{H}_{p_{0}}$ is Hilbert--Schmidt and for all $%
\varepsilon >0$ small enough $(\varepsilon \leq (2^{q_{0}}e\left\|
i_{p,p_{0}}\right\| _{HS})^{-1})$ there exists a constant $C_{p,\varepsilon
}>0$ with 
\begin{equation}
\left| P_{n}^{\mu }(z)\right| _{-p}\leq C_{p,\varepsilon }\,n!\,\varepsilon
^{-n}\,e\left( \varepsilon |z|_{-p}\right) ,\quad z\in \mathcal{H}_{-p,%
\QTR{mathbb}{C}}.  \label{1eq51}
\end{equation}
\end{proposition}

The following lemma describes the set of polynomials $\mathcal{P}(\mathcal{N}%
^{\prime }).$

\begin{lemma}
\label{1eq52}For any $\varphi \in \mathcal{P}(\mathcal{N}%
^{\prime })$ there exists a unique representation 
\begin{equation}
\varphi \left( x\right) =\sum_{n=0}^{N}\left\langle P_{n}^{\mu }\left(
x\right) ,\varphi ^{\left( n\right) }\right\rangle \,,\quad \varphi ^{\left(
n\right) }\in \mathcal{N}_{\QTR{mathbb}{C}}^{\widehat{\otimes }n}
\label{1eq53}
\end{equation}
and vice versa, any functional of the form (\ref{1eq53}) is a smooth
polynomial.
\end{lemma}

\subsection{\textrm{\textbf{Q}}$^{\mu }${-}system}

To give an internal description of the type (\ref{1eq53}) for $%
\mathcal{P}_{\mu }^{\prime }(\mathcal{N}^{\prime })$ we have to cons\-truct
an appropriate system of generalized functions, the \textrm{\textbf{Q}}$%
^{\mu }${-}system. We propose to construct the \textrm{\textbf{Q}}$^{\mu }${-%
}system using differential operators.

For $\Phi ^{\left( n\right) }\in \mathcal{N}_{\QTR{mathbb}{C}}^{\prime 
\widehat{\otimes }n}$ define a differential operator, $D(\Phi ^{\left(
n\right) })$, of order $n$ and constant coefficients $\Phi ^{\left( n\right)
}\in \mathcal{N}_{\QTR{mathbb}{C}}^{\prime \widehat{\otimes }n}$, such that,
applied to monomials $\langle x^{\otimes m},\varphi ^{\left( m\right)
}\rangle $, $\varphi ^{\left( m\right) }\in \mathcal{N}_{\QTR{mathbb}{C}}^{%
\widehat{\otimes }m}$,$\,m\in \QTR{mathbb}{N}$%
\begin{equation}
D\left( \Phi ^{\left( n\right) }\right) \left\langle x^{\otimes m},\varphi
^{\left( m\right) }\right\rangle =\left\{ 
\begin{array}{cc}
\displaystyle\frac{m!}{\left( m-n\right) !}\left\langle x^{\otimes \left(
m-n\right) }\widehat{\otimes }\Phi ^{\left( n\right) },\varphi ^{\left(
m\right) }\right\rangle  & \mathrm{for\ }m\geq n \\ 
0 & \mathrm{for\ }m<n
\end{array}
\right.   \label{1eq54}
\end{equation}
and extend by linearity from the monomials to $\mathcal{P}(\mathcal{N}%
^{\prime }).$

\begin{lemma}
\label{1eq55}$D(\Phi ^{\left( n\right) })$ is a continuous linear
operator from $\mathcal{P}(\mathcal{N}^{\prime })$ to $\mathcal{P}(\mathcal{N%
}^{\prime }).$
\end{lemma}

\noindent \textbf{Remark}\quad For $\Phi ^{\left( 1\right) }\in \mathcal{N}%
^{\prime }$ we have the usual G\^{a}teaux derivative as e.g., in white noise
analysis \cite{HKPS93} 
\[
D\left( \Phi ^{\left( 1\right) }\right) \varphi =D_{\Phi ^{\left( 1\right)
}}\varphi :=\frac{d}{dt}\left. \varphi \left( \cdot +t\Phi ^{\left( 1\right)
}\right) \right| _{t=0} 
\]
for $\varphi \in \mathcal{P}(\mathcal{N}^{\prime }).$ Moreover we have $%
D((\Phi ^{\left( 1\right) })^{\otimes n})=(D_{\Phi ^{\left( 1\right) }})^{n}$%
, thus $D((\Phi ^{\left( 1\right) })^{\otimes n})$ is a differential
operator of order $n.$

In view of Lemma \ref{1eq55} it is possible to define the adjoint
operator 
\[
D(\Phi ^{\left( n\right) })^{*}:\mathcal{P}_{\mu }^{\prime }(\mathcal{N}%
^{\prime })\longrightarrow \mathcal{P}_{\mu }^{\prime }(\mathcal{N}^{\prime
}),\quad \Phi ^{\left( n\right) }\in \mathcal{N}_{\QTR{mathbb}{C}}^{\prime 
\widehat{\otimes }n}. 
\]
Further we introduce the constant function ${1\in L}^{2}(\mu )\subset 
\mathcal{P}_{\mu }^{\prime }(\mathcal{N}^{\prime })$ such that ${1}\left(
x\right) \equiv 1$ for all $x\in \mathcal{N}^{\prime }$, so 
\[
\left\langle \!\left\langle 1,\varphi \right\rangle \!\right\rangle _{\mu
}=\int_{\mathcal{N}^{\prime }}\varphi \left( x\right) d\mu \left( x\right) =%
\QTR{mathbb}{E}_{\mu }\left( \varphi \right) . 
\]

Now we are ready to define the \textrm{\textbf{Q}}$^{\mu }$-system.

\begin{definition}
For any $\Phi ^{\left( n\right) }\in \mathcal{N}_{\QTR{mathbb}{C}}^{\prime 
\widehat{\otimes }n}$ we define a generalized function $Q_{n}^{\mu }(\Phi
^{\left( n\right) })\in \mathcal{P}_{\mu }^{\prime }(\mathcal{N}^{\prime })$
by 
\[
Q_{n}^{\mu }\left( \Phi ^{\left( n\right) }\right) =D\left( \Phi ^{\left(
n\right) }\right) ^{*}1{.} 
\]
\end{definition}

We want to introduce an additional formal notation which stresses the
linearity of $\Phi ^{\left( n\right) }\mapsto Q_n^\mu (\Phi ^{\left(
n\right) })\in \mathcal{P}_\mu ^{\prime }(\mathcal{N}^{\prime }):$%
\[
\left\langle Q_n^\mu ,\Phi ^{\left( n\right) }\right\rangle :=Q_n^\mu \left(
\Phi ^{\left( n\right) }\right) 
\]

\begin{example}
The simplest non trivial case can be studied using finite dimensional real
analysis. We consider the nuclear ''triple'' 
\[
\QTR{mathbb}{R}{\subseteq \QTR{mathbb}{R}\subseteq \QTR{mathbb}{R}} 
\]
where the dual pairing between a ''test function'' and a ''distribution''
degene\-rates to multiplication. On $\QTR{mathbb}{R}$ we consider a measure d%
$\mu \left( x\right) =\rho \left( x\right) dx$ where $\rho $ is a positive $%
C^{\infty }$-function on $\QTR{mathbb}{R}$ such that assumptions 1 and 2 are
fulfilled. In this setting the adjoint of the differentiation operator is
given by 
\[
\left( \frac{d}{dx}\right) ^{*}f\left( x\right) =-\left( \left( \frac{d}{dx}%
\right) +\beta \left( x\right) \right) f\left( x\right) ,\quad f\in
C^{\infty }\left( \QTR{mathbb}{R}\right) , 
\]
where $\beta $ is the logarithmic derivative of the measure $\mu $ and given
by 
\[
\beta =\frac{\rho ^{\prime }}{\rho }. 
\]
This enables us to calculate the \textrm{\textbf{Q}}$^{\mu }$-system. One
has 
\begin{eqnarray*}
Q_{n}^{\mu }\left( x\right) &=&\left( \left( \frac{d}{dx}\right) ^{*}\right)
^{n}{1} \\
\ &=&\left( -1\right) ^{n}\left( \frac{d}{dx}+\beta \left( x\right) \right)
^{n}{1} \\
\ &=&\left( -1\right) ^{n}\frac{\rho ^{\left( n\right) }\left( x\right) }{%
\rho \left( x\right) },
\end{eqnarray*}
where the last equality can be seen by simple induction (for $\rho \,$ non
smooth this construction produce \emph{\textbf{generalized}\ }functions $%
Q_{n}^{\mu }$ even in this one dimensional case).

If $\rho \left( x\right) =\frac{1}{\sqrt{2\pi }}\exp (-\frac{1}{2}x^{2})$ is
the Gaussian density, then $Q_{n}^{\mu }$ is related to the n-th Hermite
polynomial: 
\[
Q_{n}^{\mu }\left( x\right) =2^{-n/2}H_{n}\left( \frac{x}{\sqrt{2}}\right) . 
\]
\end{example}

\begin{definition}
We define the \textrm{\textbf{Q}}$^{\mu }$-system in $\mathcal{P}_{\mu
}^{\prime }(\mathcal{N}^{\prime })$ by 
\[
\mathrm{\mathbf{Q}}^{\mu }=\left\{ Q_{n}^{\mu }\left( \Phi ^{\left( n\right)
}\right) \,|\;\Phi ^{\left( n\right) }\in \mathcal{N}_{\QTR{mathbb}{C}%
}^{\prime \widehat{\otimes }n},\;n\in \QTR{mathbb}{N}_{0}\right\} , 
\]
and the pair $(\mathrm{\mathbf{P}}^{\mu },\mathrm{\mathbf{Q}}^{\mu })$ will
be called the \emph{\textbf{Appell system}} \textrm{\textbf{A}}$^{\mu }$
generated by the measure $\mu .$
\end{definition}

We have the following central property of the Appell system \textrm{\textbf{A%
}}$^{\mu }.$

\begin{theorem}
\label{1eq56}\emph{\textbf{(Biorthogonality w.r.t. $\mu $)}} 
\begin{equation}
\left\langle \!\!\left\langle Q_{n}^{\mu }(\Phi ^{(n)}),\ P_{m}^{\mu }\left(
\varphi ^{(m)}\right) \right\rangle \!\!\right\rangle _{\mu }=\delta
_{m,n}\;n!\;\langle \Phi ^{(n)},\varphi ^{(n)}\rangle  \label{1eq57}
\end{equation}
for $\Phi ^{(n)}\in \mathcal{N}_{\QTR{mathbb}{C}}^{\prime \widehat{\otimes }%
n}$ and $\varphi ^{(m)}\in \mathcal{N}_{\QTR{mathbb}{C}}^{\widehat{\otimes }%
m}$.
\end{theorem}

Now we are going to characterize the space $\mathcal{P}_\mu ^{\prime }(%
\mathcal{N}^{\prime }).$

\begin{theorem}
\label{1eq58}For all $\Phi \in \mathcal{P}_{\mu }^{\prime }(\mathcal{N}%
^{\prime })$ there exists a unique sequence $\{\Phi ^{\left( n\right)
}\,|\,n\in \QTR{mathbb}{N}_{0}\}$, $\Phi ^{\left( n\right) }\in \mathcal{N}_{%
\QTR{mathbb}{C}}^{\prime \widehat{\otimes }n}$ such that 
\begin{equation}
\Phi =\sum_{n=0}^{\infty }Q_{n}^{\mu }\left( \Phi ^{\left( n\right) }\right)
\equiv \sum_{n=0}^{\infty }\left\langle Q_{n}^{\mu },\Phi ^{\left( n\right)
}\right\rangle  \label{1eq59}
\end{equation}
and vice versa, every series of the form (\ref{1eq59}) generates a
generalized function in $\mathcal{P}_{\mu }^{\prime }(\mathcal{N}^{\prime
})$.
\end{theorem}

The proofs of this result can be found in \cite{KSWY95}.

\section{The triple $(\mathcal{N})^1\subset L^2(\mu )\subset (\mathcal{N}%
)_\mu ^{-1}$}

\subsection{Test functions}

\noindent On the space $\mathcal{P}(\mathcal{N}^{\prime })$ we can define a
system of norms using the Appell decomposition from Lemma \ref
{1eq52}. Let 
\[
\varphi \left( x\right) =\sum_{n=0}^{N}\left\langle P_{n}^{\mu }\left(
x\right) ,\varphi ^{\left( n\right) }\right\rangle \in \mathcal{P}(\mathcal{N%
}^{\prime })
\]
be given, then $\varphi ^{\left( n\right) }\in \mathcal{H}_{p,\QTR{mathbb}{C}%
}^{\widehat{\otimes }n}$ for each $p\geq 0$ ($n\in \QTR{mathbb}{N}_{0}{)}${.}
Thus we may define for any $p,q\in \QTR{mathbb}{N}$ a Hilbert norm on $%
\mathcal{P}(\mathcal{N}^{\prime })$ by 
\[
\left\| \varphi \right\| _{p,q,\mu }^{2}=\sum_{n=0}^{\infty }\left(
n!\right) ^{2}2^{nq}\left| \varphi ^{\left( n\right) }\right| _{p}^{2}
\]
The completion of $\mathcal{P}(\mathcal{N}^{\prime })$ w.r.t. $\left\| \cdot
\right\| _{p,q,\mu }$ is denoted by $(\mathcal{H}_{p})_{q,\mu }^{1}\,$.

\begin{definition}
We define 
\[
\left( \mathcal{N}\right) _{\mu }^{1}:=\,\stackunder{p,q\in \QTR{mathbb}{N}}{%
\mathrm{pr\,\lim }}(\mathcal{H}_{p})_{q,\mu }^{1}
\]
\end{definition}

This space have the following properties (for the proofs see \cite{KSWY95}
and references therein).

\begin{theorem}
$(\mathcal{N})_{\mu }^{1}$ is a nuclear space. The topology $(\mathcal{N}%
)_{\mu }^{1}$ is uniquely defined by the topology on $\mathcal{N}$: It does
not depend on the choice of the family of norms $\{\left| \cdot \right|
_{p}\}$.
\end{theorem}

\begin{theorem}
There exists $p^{\prime },q^{\prime }>0$ such that for all $p\geq p^{\prime }
$,$\;q\geq q^{\prime }$ the topological embedding $(\mathcal{H}_{p})_{q,\mu
}^{1}\subset L^{2}(\mu )$ holds.
\end{theorem}

\begin{corollary}
$(\mathcal{N})_{\mu }^{1}$ is continuously and densely embedded in $%
L^{2}(\mu )$.
\end{corollary}

\begin{theorem}
Any test function $\varphi $ in $(\mathcal{N})_{\mu }^{1}$ has a uniquely
defined extension to $\mathcal{N}_{\QTR{mathbb}{C}}^{\prime }$ as an element
of $\mathcal{E}_{\min }^{1}(\mathcal{N}_{\QTR{mathbb}{C}}^{\prime })$.
\end{theorem}

\noindent In this construction one unexpected moment was the following:

\begin{theorem}
\label{1eq1}For all measures $\mu \in \mathcal{M}_{a}(\mathcal{N}%
^{\prime })$ we have the topological identity 
\[
\left( \mathcal{N}\right) _{\mu }^{1}=\mathcal{E}_{\min }^{1}\left( \mathcal{%
N}^{\prime }\right) .
\]
\end{theorem}

Since this last theorem states that the space of test functions $(\mathcal{N}%
)_{\mu }^{1}$ is isomorphic to $\mathcal{E}_{\min }^{1}(\mathcal{N}^{\prime
})$ for all measures $\mu \in \mathcal{M}_{a}(\mathcal{N}^{\prime })$, we
will drop the subscript $\mu $. The test function space $(\mathcal{N})^{1}$
is the same for all measures $\mu \in \mathcal{M}_{a}(\mathcal{N}^{\prime })$%
.

\subsection{Distributions}

\noindent The space $(\mathcal{N})_{\mu }^{-1}$ of distributions
corresponding to the space of test functions $(\mathcal{N})^{1}$ can be
viewed as a subspace of $\mathcal{P}_{\mu }^{\prime }(\mathcal{N}^{\prime })$%
, since $\mathcal{P}(\mathcal{N}^{\prime })\subset (\mathcal{N})^{1}$
topologically, i.e., 
\[
(\mathcal{N})_{\mu }^{-1}\subset \mathcal{P}_{\mu }^{\prime }(\mathcal{N}%
^{\prime })
\]
Let us now introduce the Hilbert subspace $(\mathcal{H}_{-p})_{-q,\mu }^{-1}$
of $\mathcal{P}_{\mu }^{\prime }(\mathcal{N}^{\prime })$ for which the norm 
\[
\left\| \Phi \right\| _{-p,-q,\mu }^{2}:=\sum_{n=0}^{\infty }2^{-qn}\left|
\Phi ^{\left( n\right) }\right| _{-p}^{2}
\]
is finite. Here we used the canonical representation 
\[
\Phi =\sum_{n=0}^{\infty }Q_{n}^{\mu }\left( \Phi ^{\left( n\right) }\right)
\in \mathcal{P}_{\mu }^{\prime }(\mathcal{N}^{\prime })
\]
from Theorem \ref{1eq58}. The space $(\mathcal{H}_{-p})_{-q,\mu }^{-1}$
is the dual space of $(\mathcal{H}_{p})_{q,\mu }^{1}$ with respect to $%
L^{2}(\mu )$ (because of the biorthogonality of \textrm{\textbf{P}}$^{\mu }${%
- }and \textrm{\textbf{Q}}$^{\mu }${-}systems). By the general duality
theory 
\[
(\mathcal{N})_{\mu }^{-1}=\bigcup_{p,q\in \QTR{mathbb}{N}}\left( \mathcal{H}%
_{-p}\right) _{-q,\mu }^{-1}
\]
is the dual space of $(\mathcal{N})^{1}$ with respect to $L^{2}(\mu )$. So,
we have the topological nuclear triple 
\[
(\mathcal{N})^{1}\subset L^{2}(\mu )\subset \left( \mathcal{N}\right) _{\mu
}^{-1}.
\]
The action of 
\[
\Phi =\sum_{n=0}^{\infty }Q_{n}^{\mu }\left( \Phi ^{\left( n\right) }\right)
\in \left( \mathcal{N}\right) _{\mu }^{-1}
\]
on a test function 
\[
\varphi =\sum_{n=0}^{\infty }\left\langle P_{n}^{\mu },\varphi ^{\left(
n\right) }\right\rangle \in \left( \mathcal{N}\right) ^{1}
\]
is given by 
\[
\left\langle \!\left\langle \Phi ,\varphi \right\rangle \!\right\rangle
_{\mu }=\sum_{n=0}^{\infty }n!\left\langle \Phi ^{\left( n\right) },\varphi
^{\left( n\right) }\right\rangle .
\]

\begin{example}
\emph{\textbf{(Generalized Radon-Nikodym derivative)}}\textbf{\ }We want to
define a generalized function $\rho _{\mu }\left( z,\cdot \right) \in (%
\mathcal{N})_{\mu }^{-1},$ $z\in \mathcal{N}_{\QTR{mathbb}{C}}^{\prime }$
with the following property 
\[
\left\langle \!\left\langle \rho _{\mu }\left( z,\cdot \right) ,\varphi
\right\rangle \!\right\rangle _{\mu }=\int_{\mathcal{N}^{\prime }}\varphi
\left( x-z\right) d\mu \left( x\right) ,\;\varphi \in \left( \mathcal{N}%
\right) ^{1}.
\]
That means we have to establish the continuity of $\rho _{\mu }\left(
z,\cdot \right) $. Let $z\in \mathcal{H}_{-p,\QTR{mathbb}{C}}$. If $%
p^{\prime }\geq p$ is sufficiently large and $\epsilon >0$ is small enough,
there exists $q\in \QTR{mathbb}{N}$ and $C>0$ such that 
\begin{eqnarray*}
\left| \int_{\mathcal{N}^{\prime }}\varphi \left( x-z\right) d\mu \left(
x\right) \right|  &\leq &C\left\| \varphi \right\| _{p^{\prime },q,\mu
}\int_{\mathcal{N}^{\prime }}\exp \left( \epsilon \left| x-z\right|
_{-p^{\prime }}\right) d\mu \left( x\right)  \\
&\leq &C\left\| \varphi \right\| _{p^{\prime },q,\mu }\exp \left( \epsilon
\left| z\right| _{-p^{\prime }}\right) \int_{\mathcal{N}^{\prime }}\exp
\left( \epsilon \left| x\right| _{-p^{\prime }}\right) d\mu \left( x\right) .
\end{eqnarray*}
If $\epsilon $ is chosen sufficiently small the last integral exists. Thus
we have in fact $\rho _{\mu }\left( z,\cdot \right) \in (\mathcal{N})_{\mu
}^{-1}$. It is clear that whenever the Radon-Nikodym derivative $\frac{d\mu
(x+\xi )}{d\mu (x)}$ exists (e.g., $\xi \in \mathcal{N}$ in case $\mu $ is $%
\mathcal{N}$-quasi-invariant) it coincides with $\rho _{\mu }\left( z,\cdot
\right) $ defined above. We will show that in $(\mathcal{N})_{\mu }^{-1}$ we
have the canonical expansion 
\[
\rho _{\mu }\left( z,\cdot \right) =\sum_{n=0}^{\infty }\frac{\left(
-1\right) ^{n}}{n!}Q_{n}^{\mu }\left( z^{\otimes n}\right) .
\]
Since both sides are in $(\mathcal{N})_{\mu }^{-1}$ it is sufficient to
compare their action on a total set from $(\mathcal{N})^{1}$. For $\varphi
^{\left( n\right) }\in \mathcal{N}_{\QTR{mathbb}{C}}^{\widehat{\otimes }n}$
we have 
\begin{eqnarray*}
&&\left\langle \!\!\left\langle \rho _{\mu }\left( z,\cdot \right)
,\left\langle P_{n}^{\mu },\varphi ^{\left( n\right) }\right\rangle
\right\rangle \!\!\right\rangle _{\mu } \\
&=&\int_{\mathcal{N}^{\prime }}\left\langle P_{n}^{\mu }\left( x-z\right)
,\varphi ^{\left( n\right) }\right\rangle d\mu \left( x\right)  \\
&=&\sum_{k=0}^{n}\binom{n}{k}\left( -1\right) ^{n-k}\int_{\mathcal{N}%
^{\prime }}\left\langle P_{k}^{\mu }\left( x\right) \widehat{\otimes }%
z^{\otimes \left( n-k\right) },\varphi ^{\left( n\right) }\right\rangle d\mu
\left( x\right)  \\
&=&\left( -1\right) ^{n}\left\langle z^{\otimes n},\varphi ^{\left( n\right)
}\right\rangle  \\
&=&\left\langle \!\!\left\langle \sum_{k=0}^{\infty }\frac{1}{k!}\left(
-1\right) ^{k}Q_{k}^{\mu }\left( z^{\otimes k}\right) ,\left\langle
P_{n}^{\mu },\varphi ^{\left( n\right) }\right\rangle \right\rangle
\!\!\right\rangle _{\mu }\,,
\end{eqnarray*}
where we have used (\ref{1eq49}), (\ref{1eq50}) and the biorthogonality of 
\textrm{\textbf{P}}$^{\mu }$- and \textrm{\textbf{Q}}$^{\mu }$-systems. 
In other words, we have proven that $\rho _{\mu }\left( -z,\cdot \right) $ 
is the generating function of the \textrm{\textbf{Q}}$^{\mu }$-system 
\[
\rho _{\mu }\left( -z,\cdot \right) =\sum_{n=0}^{\infty }\frac{1}{n!}%
Q_{n}^{\mu }\left( z^{\otimes n}\right) .
\]
\end{example}

\subsection{Integral transformations}

\subsubsection{Normalized Laplace transform\textbf{\ }$S_\mu $}

We first introduce the Laplace transform of a function $\varphi \in
L^{2}(\mu )$. The global assumption $\mu \in \mathcal{M}_{a}(\mathcal{N}%
^{\prime })$ guarantees the existence of $p_{\mu }^{\prime }\in \QTR{mathbb}{%
N}{,\epsilon }_{\mu }>0$ such that 
\[
\int_{\mathcal{N}^{\prime }}\exp \left( -\epsilon _{\mu }\left| x\right|
_{-p_{\mu }}\right) d\mu \left( x\right) <\infty 
\]
by Lemma \ref{1eq11}. Thus $\exp (\langle \cdot ,\theta \rangle )\in
L^{2}(\mu )$ if $2\left| \theta \right| _{p_{\mu }^{\prime }}<\epsilon _{\mu
}$, $\theta \in \mathcal{H}_{p_{\mu }^{\prime },\QTR{mathbb}{C}}$. Then by
Cauchy-Schwarz inequality the Laplace transform defined by 
\[
L_{\mu }\varphi \left( \theta \right) :=\int_{\mathcal{N}^{\prime }}\varphi
\left( x\right) \exp \left\langle x,\theta \right\rangle d\mu \left(
x\right) 
\]
is well defined for $\varphi \in L^{2}(\mu )$, $\theta \in \mathcal{H}%
_{p_{\mu }^{\prime },\QTR{mathbb}{C}}$. Now we are interested to extend this
integral transform from $L^{2}(\mu )$ to the space of distributions $(%
\mathcal{N})_{\mu }^{-1}$.

Since our construction of test functions and distributions spaces is closely
related to \textrm{\textbf{P}}$^{\mu }$- and \textrm{\textbf{Q}}$^{\mu }$%
-systems it is useful to introduce the so called $S_{\mu }$-transform 
\[
S_{\mu }\varphi \left( \theta \right) :=\frac{L_{\mu }\varphi \left( \theta
\right) }{l_{\mu }\left( \theta \right) }=\int_{\mathcal{N}^{\prime
}}\varphi \left( x\right) e_{\mu }\left( \theta ;x\right) d\mu \left(
x\right) .
\]

The $\mu $-exponential $e_{\mu }\left( \theta ;\cdot \right) $ is not a test
function in $(\mathcal{N})^{1}$, see \cite[Example 6]{KSWY95}, so the
definition of the $S_{\mu }$-transform of a distribution $\Phi \in (\mathcal{%
N})_{\mu }^{-1}$ must be more careful. Every such $\Phi $ is of finite
order, i.e., $\exists p,q\in \QTR{mathbb}{N}$ such that $\Phi \in (\mathcal{H%
}_{-p})_{-q,\mu }^{-1}$ and $e_{\mu }\left( \theta ;\cdot \right) $ is in
the corresponding dual space $(\mathcal{H}_{p})_{q,\mu }^{1}$ if $\theta \in 
\mathcal{H}_{p,\QTR{mathbb}{C}}$ is such that $2^{q}\left| \theta \right|
_{p}^{2}<1$. Then we can define a consistent extension of $S_{\mu }$%
-transform. 
\[
S_{\mu }\Phi \left( \theta \right) :=\left\langle \!\left\langle \Phi
,e_{\mu }\left( \theta ,\cdot \right) \right\rangle \!\right\rangle _{\mu
}\!\,
\]
if $\theta $ is chosen in the above way. The biorthogonality of \textrm{%
\textbf{P}}$^{\mu }${- and }\textrm{\textbf{Q}}$^{\mu }$-system implies 
\[
S_{\mu }\Phi \left( \theta \right) =\sum_{n=0}^{\infty }\left\langle \Phi
^{\left( n\right) },\theta ^{\otimes n}\right\rangle ,
\]
moreover $S_{\mu }\Phi \in \mathrm{Hol}_{0}(\mathcal{N}_{\QTR{mathbb}{C}})$,
see \cite[Theorem 35]{KSWY95}.

\subsubsection{Convolution\ $C_\mu $}

We define the convolution of a function $\varphi \in (\mathcal{N})^{1}$ with
the measure $\mu $ by 
\[
C_{\mu }\varphi \left( y\right) :=\int_{\mathcal{N}^{\prime }}\varphi \left(
x+y\right) d\mu \left( x\right) ,\quad y\in \mathcal{N}^{\prime }.
\]

For any $\varphi \in (\mathcal{N})^{1}$, $z\in \mathcal{N}_{\QTR{mathbb}{C}%
}^{\prime }$, the convolution has the representation 
\[
C_{\mu }\varphi \left( z\right) =\left\langle \!\left\langle \rho _{\mu
}\left( -z,\cdot \right) ,\varphi \right\rangle \!\right\rangle _{\mu }.
\]
If $\varphi \in (\mathcal{N})^{1}$ has the canonical \textrm{\textbf{P}}$%
^{\mu }${-decomposition} 
\[
\varphi =\sum_{n=0}^{\infty }\left\langle P_{n}^{\mu },\varphi ^{\left(
n\right) }\right\rangle ,
\]
then 
\[
C_{\mu }\varphi \left( z\right) =\sum_{n=0}^{\infty }\left\langle z^{\otimes
n},\varphi ^{\left( n\right) }\right\rangle .
\]

In Gaussian analysis $C_\mu $- and $S_\mu $-transform coincide. It is a
typical non-Gaussian effect that these two transformations differ from each
other.

\subsection{Characterization theorems}

\noindent Now we will characterize the spaces of test and generalized
functions by the integral transforms introduced in the previous section.

We will start to characterize the space $(\mathcal{N})^{1}$ in terms of the
convolution $C_{\mu }$.

\begin{theorem}
The convolution $C_{\mu }$ is a topological isomorphism from $(\mathcal{N}%
)^{1}$ on $\mathcal{E}_{\mathrm{\min }}^{1}(\mathcal{N}_{\QTR{mathbb}{C}%
}^{\prime })$.
\end{theorem}

\bigskip \noindent \textbf{Remark}.\quad Since we have identified $(\mathcal{%
N})^{1}$ and $\mathcal{E}_{\mathrm{\min }}^{1}(\mathcal{N}^{\prime })$ by
Theorem \ref{1eq1}, the above assertion can be restated as follows. We
have 
\[
C_{\mu }:\mathcal{E}_{\mathrm{\min }}^{1}(\mathcal{N}^{\prime })\rightarrow 
\mathcal{E}_{\mathrm{\min }}^{1}(\mathcal{N}_{\QTR{mathbb}{C}}^{\prime }),
\]
as a topological isomorphism.

The next Theorem characterizes distributions from $(\mathcal{N})_\mu ^{-1}$
in terms of $S_\mu $-transform.

\begin{theorem}
\label{1eq2}The $S_{\mu }$-transform is a topological
isomorphism from $(\mathcal{N})_{\mu }^{-1}$ on \textrm{Hol}$_{0}(\mathcal{N}%
_{\QTR{mathbb}{C}})$.
\end{theorem}

\noindent Detailed proofs of the above theorems can be founded in 
\cite[Theorems 33, 35]{KSWY95}.

\section{Generalized Appell Systems}

\subsection{Description of the \textrm{\textbf{P}}$^{\mu ,\alpha }$-system%
\label{1eq15}}

\noindent Remember that the $\mu $-exponential is the generating function of
the \textrm{\textbf{P}}$^{\mu }$-system, i.e., if $\theta \in \mathcal{U}%
_{0}\subset \mathcal{N}_{\QTR{mathbb}{C}}$ and $z\in \mathcal{N}_{%
\QTR{mathbb}{C}}^{\prime }$, then 
\[
e_{\mu }\left( \theta ,z\right) :=\frac{\exp \left\langle z,\theta
\right\rangle }{l_{\mu }\left( \theta \right) }=\sum_{n=0}^{\infty }\frac{1}{%
n!}\left\langle P_{n}^{\mu }\left( z\right) ,\theta ^{\otimes
n}\right\rangle ,\;P_{n}^{\mu }\left( z\right) \in \mathcal{N}_{\QTR{mathbb}{%
C}}^{\prime \widehat{\otimes }n}. 
\]

In view to generalize the Appell system we consider $\alpha \in \mathrm{Hol}%
_{0}(\mathcal{N}_{\QTR{mathbb}{C}},\mathcal{N}_{\QTR{mathbb}{C}})$ an
invertible function such that $\alpha (0)=0$; moreover we have the following
decomposition 
\begin{equation}
\alpha \left( \theta \right) =\sum_{n=1}^{\infty }\frac{1}{n!}\left\langle
\alpha ^{\left( n\right) }\left( 0\right) ,\theta ^{\otimes n}\right\rangle
,\quad \theta \in \mathcal{U}_{\alpha }\subset \mathcal{N}_{\QTR{mathbb}{C}}
\label{1eq16}
\end{equation}
where $\alpha ^{\left( n\right) }\left( 0\right) \in \mathcal{N}_{%
\QTR{mathbb}{C}}^{\prime \widehat{\otimes }n}\otimes \mathcal{N}_{%
\QTR{mathbb}{C}}$ since $\alpha $ is vector valued. Analogously for the
inverse function $\alpha ^{-1}=:g_{\alpha }$, we have 
\begin{equation}
g_{\alpha }\left( \theta \right) =\sum_{n=1}^{\infty }\frac{1}{n!}%
\left\langle g_{\alpha }^{\left( n\right) }\left( 0\right) ,\theta ^{\otimes
n}\right\rangle ,\;\theta \in \mathcal{V}_{\alpha }\subset \mathcal{N}_{%
\QTR{mathbb}{C}},  \label{1eq17}
\end{equation}
where $g_{\alpha }^{\left( n\right) }\left( 0\right) \in \mathcal{N}_{%
\QTR{mathbb}{C}}^{\prime \widehat{\otimes }n}\otimes \mathcal{N}_{%
\QTR{mathbb}{C}}$. Now we introduce a new normalized exponential using the
function $\alpha $, i.e., 
\[
e_{\mu }^{\alpha }(\theta ;z):=e_{\mu }(\alpha (\theta );z)=\frac{\exp
\langle z,\alpha (\theta )\rangle }{l_{\mu }(\alpha (\theta ))}\,,\;\theta
\in \mathcal{U}_{\alpha }^{\prime }\subset \mathcal{U}_{\alpha },\;z\in 
\mathcal{N}_{\QTR{mathbb}{C}}^{\prime }. 
\]
Using the same procedure as in Section \ref{1eq43} there exist $%
P_{n}^{\mu ,\alpha }(z)\in \mathcal{N}_{\QTR{mathbb}{C}}^{\prime \widehat{%
\otimes }n}$ called \textbf{generalized Appell polynomial }or $\alpha $%
\textbf{-polynomial}\textit{\ }such that 
\begin{equation}
e_{\mu }^{\alpha }(\theta ;z)=\sum_{n=0}^{\infty }\frac{1}{n!}\langle
P_{n}^{\mu ,\alpha }(z),\theta ^{\otimes n}\rangle ,\;\theta \in \mathcal{U}%
_{\alpha }^{\prime }\,,\,z\in \mathcal{N}_{\QTR{mathbb}{C}}^{\prime },
\label{1eq18}
\end{equation}
which for fixed $z\in \mathcal{N}_{\QTR{mathbb}{C}}^{\prime }$ converges
uniformly on some neighborhood of zero on $\mathcal{N}_{\QTR{mathbb}{C}}$.
Hence we have constructed the \textrm{\textbf{P}}$^{\mu ,\alpha }$-system 
\[
\mathrm{\mathbf{P}}^{\mu ,\alpha }=\left\{ \left\langle P_{n}^{\mu ,\alpha
}\left( \cdot \right) ,\varphi _{\alpha }^{\left( n\right) }\right\rangle
\,|\,\varphi _{\alpha }^{\left( n\right) }\in \mathcal{N}_{\QTR{mathbb}{C}}^{%
\widehat{\otimes }n},\;n\in \QTR{mathbb}{N}\right\} . 
\]
In this case the related moments kernels of the measure $\mu $ are
determined by 
\[
l_{\mu }^{\alpha }\left( \theta \right) :=l_{\mu }\left( \alpha \left(
\theta \right) \right) =\sum_{n=0}^{\infty }\frac{1}{n!}\left\langle
M_{n}^{\mu ,\alpha },\theta ^{\otimes n}\right\rangle ,\;\theta \in \mathcal{%
N}_{\QTR{mathbb}{C}},\;M_{n}^{\mu ,\alpha }\in \mathcal{N}^{\prime \widehat{%
\otimes }n}. 
\]

Let us collect some properties of the polynomials $P_{n}^{\mu ,\alpha }(z)$.

\begin{proposition}
\label{1eq19}For $z,w\in \mathcal{N}^{\prime }$,$\ n\in \QTR{mathbb}{N}$
the following holds\\[6mm]
(P$_{\alpha }$1)\\[-11mm]
\begin{equation}
P_{n}^{\mu ,\alpha }\left( z\right) =\sum_{m=1}^{n}\frac{1}{m!}\left\langle
P_{m}^{\mu }\left( z\right) ,A_{n}^{m}\right\rangle ,  \label{1eq20}
\end{equation}
where $A_{n}^{m}$ are related to the kernels of $\alpha $ and are given in
the proof, see (\ref{1eq28}) below;\\[5mm]
(P$_{\alpha }$2)\\[-11mm]
\begin{equation}
z^{\otimes n}=\sum_{k=0}^{n}\sum_{m=0}^{k}\binom{n}{k}\frac{1}{m!}%
\left\langle P_{m}^{\mu ,\alpha }(z),B_{k}^{m}\right\rangle \widehat{\otimes 
}M_{n-k}^{\mu },  \label{1eq21}
\end{equation}
where $B_{k}^{m}$ are related with the kernels of $g_{\alpha }$ and are
given in the proof, see (\ref{1eq29}) below;\\[5mm]
(P$_{\alpha }$3)\\[-11mm]
\begin{equation}
P_{n}^{\mu ,\alpha }\left( z+w\right) =\sum_{k+l+m=n}\frac{n!}{k!l!m!}%
P_{k}^{\mu ,\alpha }\left( z\right) \widehat{\otimes }P_{l}^{\mu ,\alpha
}\left( w\right) \widehat{\otimes }M_{m}^{\mu ,\alpha }.  \label{1eq22}
\end{equation}
\\[5mm]
(P$_{\alpha }$4)\\[-11mm]
\begin{equation}
P_{n}^{\mu ,\alpha }\left( z+w\right) =\sum_{k=0}^{n}\binom{n}{k}P_{k}^{\mu
,\alpha }\left( z\right) \widehat{\otimes }P_{n-k}^{\delta _{0},\alpha
}\left( w\right) .  \label{1eq23}
\end{equation}
\ (P$_{\alpha }$5) Further, we observe 
\begin{equation}
\QTR{mathbb}{E}_{\mu }(\langle P_{m}^{\mu ,\alpha }(\cdot ),\varphi _{\alpha
}^{(m)}\rangle )=0\quad \mathrm{for\quad }m\neq 0,\;\varphi _{\alpha
}^{(m)}\in \mathcal{N}_{\QTR{mathbb}{C}}^{\widehat{\otimes }m}.
\label{1eq24}
\end{equation}
(P$_{\alpha }$6) For all $p^{\prime }>p$ such that the embedding $\mathcal{H}%
_{p^{\prime }}\hookrightarrow \mathcal{H}_{p}$ is of Hilbert-Schmidt class
and for all $\epsilon >0$ there exist ${\sigma _{\epsilon }}>0$ such that 
\begin{equation}
\left| P_{n}^{\mu ,\alpha }(z)\right| _{-p^{\prime }}\leq 2\,n!{\sigma }%
_{\epsilon }^{-n}\,\,\exp \left( \varepsilon |z|_{-p}\right) ,\quad z\in 
\mathcal{H}_{-p^{\prime },\QTR{mathbb}{C}},n\in \QTR{mathbb}{N}_{0},
\label{1eq25}
\end{equation}
where $\sigma _{\epsilon }$ is chosen in such a way that $\left| \alpha
\left( \theta \right) \right| \leq \epsilon $ and $\left| l_{\mu }\left(
\alpha \left( \theta \right) \right) \right| \geq 1/2$ for $\left| \theta
\right| _{p}={\sigma _{\epsilon }}${.}
\end{proposition}

\noindent \textbf{Proof. }(P$_{\alpha }$1) Analogously with (\ref{1eq46})
we have 
\begin{equation}
e_{\mu }^{\alpha }\left( \theta ;z\right) :=\frac{\exp \left\langle z,\alpha
\left( \theta \right) \right\rangle }{l_{\mu }\left( \alpha \left( \theta
\right) \right) }=\sum_{m=0}^{\infty }\frac{1}{m!}\left\langle P_{m}^{\mu
}\left( z\right) ,\alpha \left( \theta \right) ^{\otimes m}\right\rangle .
\label{1eq26}
\end{equation}
Using the representation from (\ref{1eq16}) we compute $\alpha
\left( \theta \right) ^{\otimes m}:$%
\begin{eqnarray}
\alpha \left( \theta \right) ^{\otimes m} &=&\sum_{l=1}^{\infty }\frac{1}{l!}%
\left\langle \alpha ^{\left( l\right) }\left( 0\right) ,\theta ^{\otimes
l}\right\rangle \otimes \cdots \otimes \sum_{l=1}^{\infty }\frac{1}{l!}%
\left\langle \alpha ^{\left( l\right) }\left( 0\right) ,\theta ^{\otimes
l}\right\rangle  \nonumber \\
&=&\sum_{l_{1},{\ldots },l_{m}=1}^{\infty }\frac{1}{l_{1}!\cdots l_{m}!}%
\left\langle \alpha ^{\left( l_{1}\right) }\left( 0\right) \otimes \cdots
\otimes \alpha ^{\left( l_{m}\right) }\left( 0\right) ,\theta ^{\otimes
\left( l_{1}+{\ldots }+l_{m}\right) }\right\rangle  \nonumber \\
&=&\sum_{n=1}^{\infty }\frac{1}{n!}\left\langle A_{n}^{m},\theta ^{\otimes
n}\right\rangle ,  \label{1eq27}
\end{eqnarray}
where 
\begin{equation}
A_{n}^{m}=\left\{ 
\begin{array}{cc}
\displaystyle \sum_{l_{1}+{\ldots }+l_{m}=n}\frac{n!}{l_{1}!{\cdots }l_{m}!}%
\alpha ^{\left( l_{1}\right) }\left( 0\right) \otimes \cdots \otimes \alpha
^{\left( l_{m}\right) }\left( 0\right) & \mathrm{for\;}n\geq m \\ 
0 & \mathrm{for\;}n<m
\end{array}
\right. .  \label{1eq28}
\end{equation}
Now we introduce (\ref{1eq27}) in (\ref{1eq26}) to obtain 
\begin{eqnarray*}
e_{\mu }^{\alpha }\left( \theta ;z\right) &=&\sum_{m=0}^{\infty }\frac{1}{m!}%
\left\langle P_{m}^{\mu }\left( z\right) ,\sum_{n=1}^{\infty }\frac{1}{n!}%
\left\langle A_{n}^{m},\theta ^{\otimes n}\right\rangle \right\rangle \\
\ &=&\sum_{n=1}^{\infty }\frac{1}{n!}\left\langle \sum_{m=0}^{n}\frac{1}{m!}%
\left\langle P_{m}^{\mu }\left( z\right) ,A_{n}^{m}\right\rangle ,\theta
^{\otimes n}\right\rangle .
\end{eqnarray*}
By definition 
\[
e_{\mu }^{\alpha }\left( \theta ;z\right) =\sum_{n=0}^{\infty }\frac{1}{n!}%
\left\langle P_{n}^{\mu ,\alpha }\left( z\right) ,\theta ^{\otimes
n}\right\rangle , 
\]
so we conclude that 
\[
P_{n}^{\mu ,\alpha }\left( z\right) =\sum_{m=1}^{n}\frac{1}{m!}\left\langle
P_{m}^{\mu }\left( z\right) ,A_{n}^{m}\right\rangle . 
\]

\noindent (P$_{\alpha }$2) Since $\theta =\alpha (g_{\alpha }\left( \theta
\right) )$ we have 
\[
e_{\mu }\left( \theta ,z\right) =\sum_{n=0}^{\infty }\frac{1}{n!}%
\left\langle P_{n}^{\mu ,\alpha }\left( z\right) ,g_{\alpha }\left( \theta
\right) ^{\otimes n}\right\rangle . 
\]
Having in mind (\ref{1eq17}) we first compute $g_{\alpha }\left(
\theta \right) ^{\otimes n}$: 
\begin{eqnarray*}
g_{\alpha }\left( \theta \right) ^{\otimes n} &=&\sum_{l=1}^{\infty }\frac{1%
}{l!}\left\langle g_{\alpha }^{\left( l\right) }\left( 0\right) ,\theta
^{\otimes l}\right\rangle \otimes {\cdots }\otimes \sum_{l=1}^{\infty }\frac{%
1}{l!}\left\langle g_{\alpha }^{\left( l\right) }\left( 0\right) ,\theta
^{\otimes l}\right\rangle \\
\ &=&\sum_{l_{1},{\ldots },l_{n}=1}^{\infty }\frac{1}{l_{1}!{\cdots }l_{n}!}%
\left\langle g_{\alpha }^{\left( l_{1}\right) }\left( 0\right) \otimes {%
\cdots }\otimes g_{\alpha }^{\left( l_{n}\right) }\left( 0\right) ,\theta
^{\otimes \left( l_{1}+{\ldots }+l_{n}\right) }\right\rangle \\
\ &=&\sum_{m=1}^{\infty }\frac{1}{m!}\left\langle B_{m}^{n},\theta ^{\otimes
m}\right\rangle ,
\end{eqnarray*}
where 
\begin{equation}
B_{m}^{n}=\left\{ 
\begin{array}{cc}
\displaystyle \sum_{l_{1}+{\ldots }+l_{n}=m}\frac{m!}{l_{1}!{\cdots }l_{n}!}%
g_{\alpha }^{\left( l_{1}\right) }\left( 0\right) \otimes {\cdots }\otimes
g_{\alpha }^{\left( l_{n}\right) }\left( 0\right) & \mathrm{for\;}m\geq n \\ 
0 & \mathrm{for\;}m<n
\end{array}
\right. .  \label{1eq29}
\end{equation}
Hence 
\begin{eqnarray*}
e_{\mu }\left( \theta ,z\right) &=&\sum_{n=0}^{\infty }\frac{1}{n!}%
\left\langle P_{n}^{\mu ,\alpha }\left( z\right) ,\sum_{m=1}^{\infty }\frac{1%
}{m!}\left\langle B_{m}^{n},\theta ^{\otimes m}\right\rangle \right\rangle \\
\ &=&\sum_{m=1}^{\infty }\frac{1}{m!}\left\langle \sum_{n=0}^{m}\frac{1}{n!}%
\left\langle P_{n}^{\mu ,\alpha }\left( z\right) ,B_{m}^{n}\right\rangle
,\theta ^{\otimes m}\right\rangle .
\end{eqnarray*}
On the other hand 
\[
e_{\mu }\left( \theta ,z\right) =\sum_{n=0}^{\infty }\frac{1}{n!}%
\left\langle P_{n}^{\mu }\left( z\right) ,\theta ^{\otimes n}\right\rangle , 
\]
so we conclude that 
\begin{equation}
P_{m}^{\mu }\left( z\right) =\sum_{n=1}^{m}\frac{1}{n!}\left\langle
P_{n}^{\mu ,\alpha }\left( z\right) ,B_{m}^{n}\right\rangle .  \label{1eq30}
\end{equation}
The result follows using property (P2) of the polynomials $P_{n}^{\mu
}\left( z\right) $.

\noindent (P$_{\alpha }$3) Let us start from the equation of the generating
functions 
\[
e_{\mu }^{\alpha }\left( \theta ,z+w\right) =e_{\mu }^{\alpha }\left( \theta
,z\right) e_{\mu }^{\alpha }\left( \theta ,w\right) l_{\mu }^{\alpha }\left(
\theta \right) . 
\]
This implies 
\begin{eqnarray*}
&&\sum_{n=0}^{\infty }\frac{1}{n!}\left\langle P_{n}^{\mu ,\alpha }\left(
z+w\right) ,\theta ^{\otimes n}\right\rangle \\
&=&\sum_{k,l,m=0}^{\infty }\frac{1}{k!l!m!}\left\langle P_{k}^{\mu ,\alpha
}\left( z\right) \widehat{\otimes }P_{l}^{\mu ,\alpha }\left( w\right) 
\widehat{\otimes }M_{m}^{\mu ,\alpha },\theta ^{\otimes \left( k+l+m\right)
}\right\rangle ,
\end{eqnarray*}
from this (P$_{\alpha }$3) follows immediately.

\noindent (P$_{\alpha }$4) We note that 
\[
e_{\mu }^{\alpha }\left( \theta ;z+w\right) =e_{\mu }^{\alpha }\left( \theta
;z\right) \exp \left\langle w,\alpha \left( \theta \right) \right\rangle
,\quad \theta \in \mathcal{U}_{0}\subset \mathcal{N}_{\QTR{mathbb}{C}}. 
\]
Now, since $l_{\delta _{0}}\left( \theta \right) =1$, we have the following
decomposition 
\begin{equation}
\exp \left\langle w,\alpha \left( \theta \right) \right\rangle
=\sum_{n=0}^{\infty }\frac{1}{n!}\left\langle P_{n}^{\delta _{0},\alpha
}\left( w\right) ,\theta ^{\otimes n}\right\rangle ,  \label{1eq31}
\end{equation}
where for $\alpha \equiv $ id, $P_{n}^{\delta _{0},\alpha }\left( w\right)
=w^{\otimes n}$. The result follows as done in (P$_{\alpha }$3).

\noindent (P$_{\alpha }$5) To see this we use, $\theta \in \mathcal{N}_{%
\QTR{mathbb}{C}}$, 
\[
\sum_{n=0}^{\infty }\frac{1}{n!}\QTR{mathbb}{E}_{\mu }\left( \left\langle
P_{m}^{\mu ,\alpha }\left( \cdot \right) ,\theta ^{\otimes n}\right\rangle
\right) =\QTR{mathbb}{E}_{\mu }\left( e_{\mu }^{\alpha }\left( \theta ;\cdot
\right) \right) =\frac{\QTR{mathbb}{E}_{\mu }\left( \exp \left\langle \cdot
,\alpha \left( \theta \right) \right\rangle \right) }{l_{\mu }\left( \alpha
\left( \theta \right) \right) }=1. 
\]
Then the polarization identity and a comparison of coefficients give the
result.

\noindent (P$_{\alpha }$6) Using the definition of $P_{n}^{\mu ,\alpha }$
and Cauchy's inequality for Taylor series we have 
\begin{eqnarray*}
\left| \left\langle P_{n}^{\mu ,\alpha }\left( z\right) ,\theta ^{\otimes
n}\right\rangle \right| &=&n!\left| \widehat{d^{n}e_{\mu }^{\alpha }\left(
0;z\right) }\left( \theta \right) \right| _{-p} \\
\ &\leq &n!\frac{1}{{\sigma }_{\epsilon }^{n}}\sup_{\left| \theta \right|
_{p}={\sigma }_{\epsilon }}\frac{\exp \left( \left| \alpha \left( \theta
\right) \right| _{p}\left| z\right| _{-p}\right) }{\left| l_{\mu }\left(
\alpha \left( \theta \right) \right) \right| }\left| \theta \right| _{p}^{n}
\\
\ &\leq &2n!{\sigma }_{\epsilon }^{-n}\exp \left( \epsilon \left| z\right|
_{-p}\right) \left| \theta \right| _{p}^{n}.
\end{eqnarray*}
The result follows by polarization and kernel theorem.\hfill $\blacksquare $%
\bigskip

Let us give a concrete example which furnish good arguments to use the 
\textrm{\textbf{P}}$^{\mu ,\alpha }$-system.

\begin{example}
\label{1eq32}\emph{\textbf{(Poisson noise)} }Let us consider the classical
(real) Schwartz triple 
\[
S\left( \QTR{mathbb}{R}\right) \subset L^{2}\left( \QTR{mathbb}{R}\right)
\subset S^{\prime }\left( \QTR{mathbb}{R}\right) . 
\]
The \emph{\textbf{Poisson white noise measure}} $\pi $ is defined as a
probability measure on $\mathcal{C}_{\sigma }(S^{\prime }(\QTR{mathbb}{R}{))}
$ with Laplace transform 
\[
l_{\pi }\left( \theta \right) =\exp \left[ \int_{\QTR{mathbb}{R}}\left( \exp
\theta \left( t\right) -1\right) dt\right] =\exp \left[ \left\langle \exp
\theta \left( \cdot \right) -1,1\right\rangle \right] ,\quad \theta \in S_{%
\QTR{mathbb}{C}}\left( \QTR{mathbb}{R}\right) , 
\]
see e.g., \cite{GV68}. It is not hard to see that $l_{\pi }$ is a
holomorphic function on $S_{\QTR{mathbb}{C}}(\QTR{mathbb}{R}{)}${,} so
assumption 1 is satisfied. But to check Assumption 2, we need additional
considerations.

First of all we remark that for any $\xi \in S(\QTR{mathbb}{R}{)}${,} $\xi
\neq 0$ the measure $\pi $ and $\pi (\cdot +\xi )$ are orthogonal (see \cite
{GGV75} for a detailed analysis). It means that $\pi $ is not $S(%
\QTR{mathbb}{R}{)}$-quasi-invariant and the note after Assumption 2 is not
applicable now.

Let some $\varphi \in \mathcal{P}(S^{\prime }(\QTR{mathbb}{R}{))}${,} $%
\varphi =0$ $\pi $-a.s. be given. We need to show that then $\varphi \equiv
0 $. To this end we will introduce a system of orthogonal polynomials in the
space $L^{2}(S^{\prime }(\QTR{mathbb}{R}{),\pi )}$ {which can be constructed
in the following way. The mapping} 
\[
\theta \left( \cdot \right) \mapsto \alpha \left( \theta \right) \left(
\cdot \right) =\log \left( 1+\theta \left( \cdot \right) \right) \in S_{%
\QTR{mathbb}{C}}\left( \QTR{mathbb}{R}\right) ,\quad \theta \in S_{%
\QTR{mathbb}{C}}\left( \QTR{mathbb}{R}\right) 
\]
is holomorphic on a neighborhood $\mathcal{U}\subset S_{\QTR{mathbb}{C}}(%
\QTR{mathbb}{R}{)}${,} $0\in \mathcal{U}$. Then 
\[
e_{\pi }^{\alpha }\left( \theta ;x\right) =\frac{\exp \left\langle x,\alpha
\left( \theta \right) \right\rangle }{l_{\pi }\left( \alpha \left( \theta
\right) \right) }=\exp \left[ \left\langle x,\alpha \left( \theta \right)
\right\rangle -\left\langle \theta ,1\right\rangle \right] ,\quad \theta \in 
\mathcal{U},\;x\in S^{\prime }\left( \QTR{mathbb}{R}\right) 
\]
is a holomorphic function on $\mathcal{U}$ for any $x\in S^{\prime }(%
\QTR{mathbb}{R}{)}${.} The Taylor decomposition and the kernel theorem give 
\[
e_{\pi }^{\alpha }\left( \theta ;x\right) =\sum_{n=0}^{\infty }\frac{1}{n!}%
\left\langle C_{n}\left( x\right) ,\theta ^{\otimes n}\right\rangle , 
\]
where $C_{n}:S^{\prime }(\QTR{mathbb}{R}{)\rightarrow }S^{\prime }(%
\QTR{mathbb}{R}{)}^{\widehat{\otimes }n}$ are polynomial mappings. For $%
\varphi ^{\left( n\right) }\in S_{\QTR{mathbb}{C}}(\QTR{mathbb}{R}{)}^{%
\widehat{\otimes }n}$, $n\in \QTR{mathbb}{N}_{0}$, we define \textit{%
Charlier polynomials} 
\[
x\mapsto C_{n}\left( \varphi ^{\left( n\right) };x\right) :=\left\langle
C_{n}\left( x\right) ,\varphi ^{\left( n\right) }\right\rangle \in 
\QTR{mathbb}{C}{,\;x\in }S^{\prime }\left( \QTR{mathbb}{R}\right) . 
\]
Due to \cite{I88}, \cite{IK88} we have the following orthogonality property: 
\[
\forall \varphi ^{\left( n\right) }\in S_{\QTR{mathbb}{C}}\left( 
\QTR{mathbb}{R}\right) ^{\widehat{\otimes }n},\;\forall \psi ^{\left(
m\right) }\in S_{\QTR{mathbb}{C}}\left( \QTR{mathbb}{R}\right) ^{\widehat{%
\otimes }m} 
\]
\[
\int C_{n}\left( \varphi ^{\left( n\right) }\right) C_{m}\left( \psi
^{\left( m\right) }\right) d\pi =\delta _{nm}n!\left\langle \varphi ^{\left(
n\right) },\psi ^{\left( n\right) }\right\rangle . 
\]
Now the rest is simple. Any continuous polynomial $\varphi $ has a uniquely
defined decomposition 
\[
\varphi \left( x\right) =\sum_{n=0}^{N}\left\langle C_{n}\left( x\right)
,\varphi ^{\left( n\right) }\right\rangle ,\quad x\in S^{\prime }\left( 
\QTR{mathbb}{R}\right) , 
\]
where $\varphi ^{\left( n\right) }\in S_{\QTR{mathbb}{C}}(\QTR{mathbb}{R}{)}%
^{\widehat{\otimes }n}$. If $\varphi =0$ $\pi $-a.e., then 
\[
\left\| \varphi \right\| _{L^{2}(\pi )}^{2}=\sum_{n=0}^{N}n!\left\langle
\varphi ^{\left( n\right) },\overline{\varphi ^{\left( n\right) }}%
\right\rangle =0. 
\]
Hence $\varphi ^{\left( n\right) }=0$, $n=0,{\ldots },N$, i.e., $\varphi
\equiv 0$. So Assumption 2 is satisfied.
\end{example}

\begin{lemma}
For any $\varphi \in \mathcal{P}(\mathcal{N}^{\prime })$ there exists a
unique representation 
\begin{equation}
\varphi \left( x\right) =\sum_{n=0}^{N}\left\langle P_{n}^{\mu ,\alpha
}\left( x\right) ,\varphi _{\alpha }^{\left( n\right) }\right\rangle
\,,\quad \varphi _{\alpha }^{\left( n\right) }\in \mathcal{N}_{\QTR{mathbb}{C%
}}^{\widehat{\otimes }n}  \label{1eq33}
\end{equation}
and vice versa, any functional of the form (\ref{1eq33}) is a
smooth polynomial.
\end{lemma}

\noindent \textbf{Proof. }The representation from Definition \ref{1eq14} and
equation (\ref{1eq33}) can be transformed into one another using (%
\ref{1eq20}) and (\ref{1eq21}).\hfill $\blacksquare $

\subsection{Description of the $\mathrm{\mathbf{Q}}^{\mu ,\alpha }$-system}

\subsubsection{Using $S_\mu $-transform}

\noindent By assumption we know that $\alpha $ is invertible with inverse
given by $g_{\alpha }$ and $\alpha \left( \theta \right) \in \mathcal{V}%
_{\alpha }\subset \mathcal{N}_{\QTR{mathbb}{C}}$, $\forall \theta \in 
\mathcal{U}_{\alpha }$. For given $\Phi _{\alpha }^{\left( n\right) }\in 
\mathcal{N}_{\QTR{mathbb}{C}}^{\prime \widehat{\otimes }n}$ we define a
generalized function $Q_{n}^{\mu ,\alpha }(\Phi _{\alpha }^{\left( n\right)
})$ via the $S_{\mu }$-transform 
\begin{equation}
S_{\mu }\left( Q_{n}^{\mu ,\alpha }\left( \Phi _{\alpha }^{\left( n\right)
}\right) \right) \left( \theta \right) :=\left\langle \Phi _{\alpha
}^{\left( n\right) },g_{\alpha }\left( \theta \right) ^{\otimes
n}\right\rangle ,\quad \,\theta \in \mathcal{V}_{\alpha }.
\label{1eq34}
\end{equation}

\subsubsection{Using differential operators}

\noindent Using the kernels $g_{\alpha }^{\left( n\right) }\left( 0\right)
\, $ of $g_{\alpha }$, see (\ref{1eq17}), we define a \textit{%
differential operator (of infinite order)} from $\mathcal{P}(\mathcal{N}%
^{\prime })$ to $\mathcal{P}(\mathcal{N}^{\prime })\otimes \mathcal{N}_{%
\QTR{mathbb}{C}}$ as follows 
\[
G_{\alpha }=\sum_{n=0}^{\infty }\frac{1}{n!}\left\langle g_{\alpha }^{\left(
n\right) }\left( 0\right) ,\bigtriangledown ^{\otimes n}\right\rangle , 
\]
such that, if $\varphi \in \mathcal{P}(\mathcal{N}^{\prime })$ and $\xi \in 
\mathcal{N}_{\QTR{mathbb}{C}}^{\prime }$ we have 
\[
G_{\alpha }^{\xi }\left( \varphi \right) \left( x\right) :=\left\langle \xi
,G_{\alpha }(\varphi )\left( x\right) \right\rangle =\sum_{n=0}^{\infty }%
\frac{1}{n!}\left\langle \xi ,\left\langle g_{\alpha }^{\left( n\right)
}\left( 0\right) ,\bigtriangledown ^{\otimes n}\varphi \left( x\right)
\right\rangle \right\rangle ,\quad x\in \mathcal{N}^{\prime }, 
\]
i.e., $G_{\alpha }^{\xi }:\mathcal{P}(\mathcal{N}^{\prime })\rightarrow 
\mathcal{P}(\mathcal{N}^{\prime })$ and formally $G_{\alpha }:=g_{\alpha
}\left( \bigtriangledown \right) $.

Let we state the following useful lemma.

\begin{lemma}
For all $\xi \in \mathcal{N}_{\QTR{mathbb}{C}}^{\prime }$, $x\in \mathcal{N}%
^{\prime }$ and $\theta \in \mathcal{N}_{\QTR{mathbb}{C}}$ we have 
\[
\left\langle \xi ,g_{\alpha }\left( \bigtriangledown \right) \right\rangle
\left( \exp \left\langle x,\theta \right\rangle \right) =\left\langle \xi
,g_{\alpha }\left( \theta \right) \right\rangle \exp \left\langle x,\theta
\right\rangle .
\]
\end{lemma}

\noindent \textbf{Proof. }Using the representation given in (\ref
{1eq17}) we have 
\[
\left\langle \xi ,g_{\alpha }\left( \bigtriangledown \right) \right\rangle
=\sum_{n=0}^{\infty }\frac{1}{n!}\left\langle g_{\alpha ,\xi }^{\left(
n\right) }\left( 0\right) ,\nabla ^{\otimes n}\right\rangle ,\quad g_{\alpha
,\xi }^{\left( n\right) }\left( 0\right) =\left\langle g_{\alpha }^{\left(
n\right) }\left( 0\right) ,\xi \right\rangle \in \mathcal{N}_{\QTR{mathbb}{C}%
}^{\prime \widehat{\otimes }n}.
\]
For simplicity we put $g_{\alpha ,\xi }^{\left( n\right) }\left( 0\right)
\equiv \Psi ^{\left( n\right) }$. At first we apply the operator to some
monomial. For given $\theta \in \mathcal{N}_{\QTR{mathbb}{C}},$ $m\geq n$%
\begin{eqnarray*}
\left\langle \Psi ^{\left( n\right) },\nabla ^{\otimes n}\right\rangle
\left\langle x,\theta \right\rangle ^{m} &=&\left\langle \Psi ^{\left(
n\right) },\nabla ^{\otimes n}\right\rangle \left\langle x^{\otimes
m},\theta ^{\otimes m}\right\rangle  \\
&=&m\left( m-1\right) \cdots \left( m-n+1\right) \left\langle \Psi ^{\left(
n\right) }\widehat{\otimes }x^{\otimes \left( m-n\right) },\theta ^{\otimes
m}\right\rangle  \\
&=&m\left( m-1\right) \cdots \left( m-n+1\right) \left\langle x,\theta
\right\rangle ^{m-n}\left\langle \Psi ^{\left( n\right) },\theta ^{\otimes
n}\right\rangle ,
\end{eqnarray*}
where we used (\ref{1eq54}) in the second equality. Now expand the
given function, $\exp \left\langle x,\theta \right\rangle $, in the Taylor
series and applying the above result we get 
\begin{eqnarray*}
&&\left\langle \Psi ^{\left( n\right) },\nabla ^{\otimes n}\right\rangle
\exp \left\langle x,\theta \right\rangle  \\
&=&\left\langle \Psi ^{\left( n\right) },\nabla ^{\otimes n}\right\rangle
\sum_{m=0}^{\infty }\frac{\left\langle x,\theta \right\rangle ^{m}}{m!} \\
&=&\sum_{m=n}^{\infty }\frac{m\left( m-1\right) \cdots \left( m-n+1\right) }{%
m!}\left\langle \Psi ^{\left( n\right) }\widehat{\otimes }x^{\otimes \left(
m-n\right) },\theta ^{\otimes m}\right\rangle  \\
&=&\left\langle \Psi ^{\left( n\right) },\theta ^{\otimes n}\right\rangle
\sum_{m=n}^{\infty }\frac{1}{\left( m-n\right) !}\left\langle x,\theta
\right\rangle ^{m-n} \\
&=&\left\langle \Psi ^{\left( n\right) },\theta ^{\otimes n}\right\rangle
\exp \left\langle x,\theta \right\rangle .
\end{eqnarray*}
Therefore 
\begin{eqnarray*}
\left\langle \xi ,g_{\alpha }\left( \bigtriangledown \right) \right\rangle
\left( \exp \left\langle x,\theta \right\rangle \right) 
&=&\sum_{n=0}^{\infty }\frac{1}{n!}\left\langle \Psi ^{\left( n\right)
},\nabla ^{\otimes n}\right\rangle \exp \left\langle x,\theta \right\rangle 
\\
&=&\sum_{n=0}^{\infty }\frac{1}{n!}\left\langle \Psi ^{\left( n\right)
},\theta ^{\otimes n}\right\rangle \exp \left\langle x,\theta \right\rangle 
\\
&=&\left\langle \xi ,g_{\alpha }\left( \theta \right) \right\rangle \left(
\exp \left\langle x,\theta \right\rangle \right) .
\end{eqnarray*}
\hfill $\blacksquare $

\begin{theorem}
Under the above conditions the $Q_{n}^{\mu ,\alpha }(\xi ^{\otimes n})$ are
given by 
\begin{equation}
Q_{n}^{\mu ,\alpha }\left( \xi ^{\otimes n}\right) \left( \cdot \right)
=\left( \left\langle \xi ,g_{\alpha }\left( \bigtriangledown \right)
\right\rangle ^{*n}{1}\right) \left( \cdot \right) .  \label{1eq35}
\end{equation}
\end{theorem}

\noindent \textbf{Proof. }Applying the $S_{\mu }$-transform to the r.h.s of (%
\ref{1eq35}) we have 
\begin{eqnarray}
S_{\mu }\left( \left\langle \xi ,g_{\alpha }\left( \bigtriangledown \right)
\right\rangle ^{*n}{1}\right) \left( \theta \right)  &=&\left\langle
\!\left\langle \left\langle \xi ,g_{\alpha }\left( \bigtriangledown \right)
\right\rangle ^{*n}{1},{e}_{\mu }\left( \theta ,\cdot \right) \right\rangle
\!\right\rangle _{\mu }  \nonumber \\
&=&\left\langle \!\left\langle {1},\left\langle \xi ,g_{\alpha }\left(
\bigtriangledown \right) \right\rangle ^{n}{e}_{\mu }\left( \theta ,\cdot
\right) \right\rangle \!\right\rangle _{\mu }  \nonumber \\
&=&\frac{1}{l_{\mu }\left( \theta \right) }\int_{\mathcal{N}^{\prime
}}\left\langle \xi ,g_{\alpha }\left( \bigtriangledown \right) \right\rangle
^{n}\exp \left\langle x,\theta \right\rangle \,\mathrm{d}\mu \left( x\right) 
\nonumber \\
&=&\frac{\left\langle \xi ,g_{\alpha }\left( \theta \right) \right\rangle
^{n}}{l_{\mu }\left( \theta \right) }\int_{\mathcal{N}^{\prime }}\exp
\left\langle x,\theta \right\rangle \,\mathrm{d}\mu \left( x\right)  
\nonumber \\
&=&\left\langle \xi ,g_{\alpha }\left( \theta \right) \right\rangle ^{n}.
\label{1eq36}
\end{eqnarray}
On the other hand the $S_{\mu }$-transform of the l.h.s. (\ref{1eq35}), by (%
\ref{1eq34}), is the same as (\ref{1eq36}) which prove the
result.\hfill $\blacksquare $

\begin{example}
As an illustration of $G_{\alpha }$ we use again the Poisson measure $\pi $
(see Example \ref{1eq32}) and $\alpha \left( \theta \right) \left( \cdot
\right) =\log \left( 1+\theta \left( \cdot \right) \right) $, $\theta \in S(%
\QTR{mathbb}{R}{)}$. For this choice we have 
\[
g_{\alpha }\left( \theta \right) \left( \cdot \right) =\exp \theta \left(
\cdot \right) -1=\sum_{n=1}^{\infty }\frac{\theta ^{n}\left( \cdot \right) }{%
n!}.
\]
On the other hand, from (\ref{1eq17}) we have 
\[
g_{\alpha }\left( \theta \right) \left( \cdot \right) =\sum_{n=1}^{\infty }%
\frac{1}{n!}\left\langle g_{\alpha }^{\left( n\right) }\left( 0\right)
,\theta ^{\otimes n}\right\rangle \left( \cdot \right) ,
\]
so we conclude that 
\[
g_{\alpha }^{\left( n\right) }\left( 0\right) =\delta \left( t_{1}-t\right)
\cdots \delta \left( t_{n}-t\right) .
\]
We introduce the notation of functional derivative (see \cite{IK88}), 
\[
\nabla _{\delta _{t}}\left( \theta \right) =\frac{\delta }{\delta \theta
\left( t\right) },\quad \theta \in S\left( \QTR{mathbb}{R}\right) {,\;t\in }%
\QTR{mathbb}{R}{.}
\]
With this, we easily see that for $\nabla _{h}=\left\langle \nabla
,h\right\rangle $ we have 
\[
\left( \exp \left( \nabla _{h}\right) f\right) \left( \cdot \right) =f\left(
\cdot +h\right) ,\quad f\in \mathcal{P}\left( S^{\prime }\left( \QTR{mathbb}{%
R}\right) \right) ,\;h\in S\left( \QTR{mathbb}{R}\right) .
\]
Hence 
\[
\left( g_{\alpha }\left( \nabla _{\delta _{t}}\right) \left( \theta \right)
\right) \left( f\left( \cdot \right) \right) =\left( \exp \left( \frac{%
\delta }{\delta \theta \left( t\right) }\right) -1\right) f\left( \cdot
\right) =f\left( \cdot +\delta _{t}\right) -f\left( \cdot \right) 
\]
and if $\xi \in S_{\QTR{mathbb}{C}}(\QTR{mathbb}{R}{)}$ we have 
\[
\left\langle g_{\alpha }\left( \nabla _{\delta _{t}}\right) ,\xi
\right\rangle f\left( \cdot \right) =\int_{\QTR{mathbb}{R}}\left[ f\left(
\cdot +\delta _{t}\right) -f\left( \cdot \right) \right] \xi \left( t\right)
dt.
\]
Therefore if $f\in \mathcal{P}(S^{\prime }(\QTR{mathbb}{R}{))}$ then 
\[
G_{\alpha }:f\left( \cdot \right) \longmapsto f\left( \cdot +\delta
_{t}\right) -f\left( \cdot \right) .
\]
This mapping can be considered as a \emph{\textbf{''gradient''}} operator on
the Poisson space $(S^{\prime }(\QTR{mathbb}{R}{),}\mathcal{B}(S^{\prime }(%
\QTR{mathbb}{R}{)),\pi })$.
\end{example}

\begin{definition}
We define the $\mathrm{\mathbf{Q}}^{\mu ,\alpha }$-system in $\mathcal{P}%
_{\mu }^{\prime }(\mathcal{N}^{\prime })$ by 
\[
\mathrm{\mathbf{Q}}^{\mu ,\alpha }=\left\{ Q_{n}^{\mu ,\alpha }\left( \Phi
_{\alpha }^{\left( n\right) }\right) \,|\,\Phi _{\alpha }^{\left( n\right)
}\in \mathcal{N}_{\QTR{mathbb}{C}}^{\prime \widehat{\otimes }n},\;n\in 
\QTR{mathbb}{N}_{0}\right\} ,
\]
and the pair $(\mathrm{\mathbf{P}}{^{\mu ,\alpha },\mathrm{\mathbf{Q}}}^{\mu
,\alpha })$ will be called the \emph{\textbf{generalized Appell system }}$%
\mathrm{\mathbf{A}}^{\mu ,\alpha }$ generated by the measure $\mu $ and
given mapping $\alpha \in \mathrm{Hol}_{0}(\mathcal{N}_{\QTR{mathbb}{C}},%
\mathcal{N}_{\QTR{mathbb}{C}})$.
\end{definition}

Now we are going to discuss the central property of the generalized Appell
system $\mathrm{\mathbf{A}}^{\mu ,\alpha }$.

\begin{theorem}
\label{1eq37}\emph{\textbf{(Biorthogonality of} }$\mathrm{\mathbf{Q}}%
^{\mu ,\alpha }$\emph{\ \textbf{and} }$\mathrm{\mathbf{P}}^{\mu ,\alpha }$%
\emph{\ \textbf{w.r.t.} }$\mu $\emph{\textbf{)} } 
\begin{equation}
\left\langle \!\!\left\langle Q_{n}^{\mu ,\alpha }\left( \Phi _{\alpha
}^{\left( n\right) }\right) ,P_{m}^{\mu ,\alpha }\left( \varphi _{\alpha
}^{\left( m\right) }\right) \right\rangle \!\!\right\rangle _{\mu }=\delta
_{nm}n!\left\langle \Phi _{\alpha }^{\left( n\right) },\varphi _{\alpha
}^{\left( n\right) }\right\rangle ,\,  \label{1eq38}
\end{equation}
for $\Phi _{\alpha }^{\left( n\right) }\in \mathcal{N}_{\QTR{mathbb}{C}%
}^{\prime \widehat{\otimes }n}$ and $\varphi _{\alpha }^{\left( m\right)
}\in \mathcal{N}_{\QTR{mathbb}{C}}^{\widehat{\otimes }m}$.
\end{theorem}

\noindent \textbf{Proof. }By definition of $S_{\mu }$ we have 
\[
S_{\mu }\left( Q_{n}^{\mu ,\alpha }\left( \Phi _{\alpha }^{\left( n\right)
}\right) \right) \left( \theta \right) :=\left\langle \!\!\left\langle
Q_{n}^{\mu ,\alpha }\left( \Phi _{\alpha }^{\left( n\right) }\right) ,e_{\mu
}(\theta ,\cdot )\right\rangle \!\!\right\rangle _{\mu }
\]
if we substitute $\theta \mapsto \alpha (\eta )$, then we obtain 
\begin{eqnarray*}
S_{\mu }\left( Q_{n}^{\mu ,\alpha }\left( \Phi _{\alpha }^{\left( n\right)
}\right) \right) \left( \alpha \left( \eta \right) \right)  &=&\left\langle
\!\!\left\langle Q_{n}^{\mu ,\alpha }\left( \Phi _{\alpha }^{\left( n\right)
}\right) ,e_{\mu }\left( \alpha \left( \eta \right) ,\cdot \right)
\right\rangle \!\!\right\rangle _{\mu } \\
&=&\sum_{m=0}^{\infty }\frac{1}{m!}\left\langle \!\!\left\langle Q_{n}^{\mu
,\alpha }\left( \Phi _{\alpha }^{\left( n\right) }\right) ,\left\langle
P_{m}^{\mu ,\alpha }\left( \cdot \right) ,\eta ^{\otimes m}\right\rangle
\right\rangle \!\!\right\rangle _{\mu }.
\end{eqnarray*}
Substituting of $\theta $ by $\alpha (\eta )$ in (\ref{1eq34})
give us 
\[
S_{\mu }\left( Q_{n}^{\mu ,\alpha }\left( \Phi _{\alpha }^{\left( n\right)
}\right) \right) \left( \alpha (\eta )\right) =\left\langle \Phi _{\alpha
}^{\left( n\right) },\eta ^{\otimes n}\right\rangle .
\]
Then a comparison of coefficients and the polarization identity give the
desired result.\hfill $\blacksquare $

\bigskip \ Now we characterize the space $\mathcal{P}_{\mu }^{\prime }(%
\mathcal{N}^{\prime })$.

\begin{theorem}
\label{1eq39}For all $\Phi \in \mathcal{P}_{\mu }^{\prime }(%
\mathcal{N}^{\prime })$ there exists a unique sequence $\{\Phi _{\alpha
}^{\left( n\right) }\,|\,n\in \QTR{mathbb}{N}_{0}\}$, $\Phi _{\alpha
}^{\left( n\right) }\in \mathcal{N}_{\QTR{mathbb}{C}}^{\prime \widehat{%
\otimes }n}$ such that 
\begin{equation}
\Phi =\sum_{n=0}^{\infty }Q_{n}^{\mu ,\alpha }\left( \Phi _{\alpha }^{\left(
n\right) }\right) \equiv \sum_{n=0}^{\infty }\left\langle Q_{n}^{\mu ,\alpha
},\Phi _{\alpha }^{\left( n\right) }\right\rangle   \label{1eq40}
\end{equation}
and vice versa, every series of the form (\ref{1eq40}) generates a
generalized function in $\mathcal{P}_{\mu }^{\prime }(\mathcal{N}^{\prime })$%
.
\end{theorem}

\noindent \textbf{Proof. }For $\Phi \in \mathcal{P}_{\mu }^{\prime }(%
\mathcal{N}^{\prime })$ we can uniquely define $\Phi _{\alpha }^{\left(
n\right) }\in \mathcal{N}_{\QTR{mathbb}{C}}^{\prime \widehat{\otimes }n}$ by 
\[
\left\langle \Phi _{\alpha }^{\left( n\right) },\varphi _{\alpha }^{\left(
n\right) }\right\rangle :=\frac{1}{n!}\left\langle \!\!\left\langle \Phi
,\left\langle P_{n}^{\mu ,\alpha },\varphi _{\alpha }^{\left( n\right)
}\right\rangle \right\rangle \!\!\right\rangle _{\mu }\,,\quad \varphi
_{\alpha }^{\left( n\right) }\in \mathcal{N}_{\QTR{mathbb}{C}}^{\widehat{%
\otimes }n}, 
\]
which is well defined since $\langle P_{n}^{\mu ,\alpha },\varphi _{\alpha
}^{\left( n\right) }\rangle \in \mathcal{P}(\mathcal{N}^{\prime })$. The
continuity of $\varphi _{\alpha }^{\left( n\right) }\mapsto \langle \Phi
_{\alpha }^{\left( n\right) },\varphi _{\alpha }^{\left( n\right) }\rangle $
follows from the continuity of $\varphi \mapsto \langle \!\langle \Phi
,\varphi \rangle \!\rangle \!_{\mu },\;\varphi \in \mathcal{P}(\mathcal{N}%
^{\prime })$. This implies that 
\[
\varphi \longmapsto \sum_{n=0}^{\infty }n!\left\langle \Phi _{\alpha
}^{\left( n\right) },\varphi _{\alpha }^{\left( n\right) }\right\rangle 
\]
is continuous on $\mathcal{P}(\mathcal{N}^{\prime })$. This defines a
generalized function in $\mathcal{P}_{\mu }^{\prime }(\mathcal{N}^{\prime })$%
, which we denote by 
\[
\sum_{n=0}^{\infty }Q_{n}^{\mu ,\alpha }\left( \Phi _{\alpha }^{\left(
n\right) }\right) . 
\]
In view of Theorem \ref{1eq37} it is easy to see that 
\[
\Phi =\sum_{n=0}^{\infty }Q_{n}^{\mu ,\alpha }\left( \Phi _{\alpha }^{\left(
n\right) }\right) . 
\]

To see the converse consider a series of the form (\ref{1eq40}) and $%
\varphi \in \mathcal{P}(\mathcal{N}^{\prime })$. Then there exists $\varphi
_{\alpha }^{\left( n\right) }\in \mathcal{N}_{\QTR{mathbb}{C}}^{\widehat{%
\otimes }n}$, $n\in \QTR{mathbb}{N}$ and $N\in \QTR{mathbb}{N}$ such that we
have the representation 
\[
\varphi =\sum_{n=0}^{N}P_{n}^{\mu ,\alpha }\left( \varphi _{\alpha }^{\left(
n\right) }\right) . 
\]
So we have 
\[
\left\langle \!\!\left\langle \sum_{n=0}^{\infty }Q_{n}^{\mu ,\alpha }\left(
\Phi _{\alpha }^{\left( n\right) }\right) ,\varphi \right\rangle
\!\!\right\rangle _{\mu }=\sum_{n=0}^{N}n!\left\langle \Phi _{\alpha
}^{\left( n\right) },\varphi _{\alpha }^{\left( n\right) }\right\rangle , 
\]
because of Theorem \ref{1eq37}. The continuity of 
\[
\varphi \longmapsto \!\!\left\langle \!\!\left\langle \sum_{n=0}^{\infty
}Q_{n}^{\mu ,\alpha }\left( \Phi _{\alpha }^{\left( n\right) }\right)
,\varphi \right\rangle \!\!\right\rangle _{\mu } 
\]
follows because $\varphi _{\alpha }^{\left( n\right) }\mapsto \langle \Phi
_{\alpha }^{\left( n\right) },\varphi _{\alpha }^{\left( n\right) }\rangle $
is continuous for all $n\in \QTR{mathbb}{N}${.}\hfill $\blacksquare $

\section{Test functions on a linear space with measure}

\subsection{Test functions spaces}

We will construct the test function space $(\mathcal{N})_{\mu ,\alpha }^{1}$
using \textrm{\textbf{P}}$^{\mu ,\alpha }$-system and study some properties.
On the space $\mathcal{P}(\mathcal{N}^{\prime })$ we can define a system of
norms using the representation from (\ref{1eq33}) 
\[
\varphi \left( \cdot \right) =\sum_{n=0}^{N}\left\langle P_{n}^{\mu ,\alpha
}\left( \cdot \right) ,\varphi _{\alpha }^{\left( n\right) }\right\rangle , 
\]
with $\varphi _{\alpha }^{\left( n\right) }\in \mathcal{H}_{p,\QTR{mathbb}{C}%
}^{\widehat{\otimes }n}$ for each $p>0$ ($n\in \QTR{mathbb}{N}$). Thus we
may define for any $p,q\in \QTR{mathbb}{N}$ a Hilbert norm on $\mathcal{P}(%
\mathcal{N}^{\prime })$ by 
\[
\left\| \varphi \right\| _{p,q,\mu ,\alpha }^{2}=\sum_{n=0}^{N}\left(
n!\right) ^{2}2^{nq}\left| \varphi _{\alpha }^{\left( n\right) }\right|
_{p}^{2}<\infty 
\]
The completion of $\mathcal{P}(\mathcal{N}^{\prime })$ w.r.t. $\left\| \cdot
\right\| _{p,q,\mu ,\alpha }^{2}$ is called $(\mathcal{H}_{p})_{q,\mu
,\alpha }^{1}\,$.

\begin{definition}
We define 
\[
\left( \mathcal{N}\right) _{\mu ,\alpha }^{1}:=\stackunder{p,q\in 
\QTR{mathbb}{N}}{\mathrm{pr\,\lim }}(\mathcal{H}_{p})_{q,\mu ,\alpha }^{1}
\]
\end{definition}

\begin{theorem}
$(\mathcal{N})_{\mu ,\alpha }^{1}$ is a nuclear space. The topology in $(%
\mathcal{N})_{\mu ,\alpha }^{1}$ is uniquely defined by the topology on $%
\mathcal{N}$. It does not depend on the choice of the family of norms $%
\{\left| \cdot \right| _{p}\}$.
\end{theorem}

\noindent \textbf{Proof. }Nuclearity of $(\mathcal{N})_{\mu ,\alpha }^{1}$
follows essentially from that of $\mathcal{N}$. For fixed $p,q$ choose $%
p^{\prime }$ such that the embedding 
\[
i_{p^{\prime },p}:\mathcal{H}_{p^{\prime }}\hookrightarrow \mathcal{H}_{p} 
\]
is Hilbert-Schmidt and consider the embedding 
\[
I_{p^{\prime },q^{\prime },p,q,\alpha }:\left( \mathcal{H}_{p^{\prime
}}\right) _{q^{\prime },\mu ,\alpha }^{1}\hookrightarrow \left( \mathcal{H}%
_{p}\right) _{q,\mu ,\alpha }^{1}. 
\]
Then $I_{p^{\prime },q^{\prime },p,q,\alpha }$ is induced by 
\[
I_{p^{\prime },q^{\prime },p,q,\alpha }\left( \varphi \right)
=\sum_{n=0}^{\infty }\left\langle P_{n}^{\mu ,\alpha },i_{p^{\prime
},p}^{\otimes n}\varphi _{\alpha }^{\left( n\right) }\right\rangle \quad 
\mathrm{for\quad }\varphi =\sum_{n=0}^{\infty }\left\langle P_{n}^{\mu
,\alpha },\varphi _{\alpha }^{\left( n\right) }\right\rangle \in (\mathcal{H}%
_{p^{\prime }})_{q^{\prime },\mu ,\alpha }^{1}. 
\]
Its Hilbert-Schmidt norm, for a given orthonormal basis of $(\mathcal{H}%
_{p^{\prime }})_{q^{\prime },\mu ,\alpha }^{1}\,$, can be estimate by 
\[
\left\| I_{p^{\prime },q^{\prime },p,q,\alpha }\right\|
_{HS}^{2}=\sum_{n=0}^{\infty }2^{n(q-q^{\prime })}\left\| i_{p^{\prime
},p}\right\| _{HS}^{2n} 
\]
which is finite for a suitably chosen $q^{\prime }$.

To prove the independence of the family of norms, let us assume that we are
given two different systems of Hilbert norms $\left| \cdot \right| _{p}$ and 
$\left| \cdot \right| _{k}^{\prime }$, such that they induce the same
topology on $\mathcal{N}$. For fixed $k$ and $l$ we have to estimate $%
\left\| \cdot \right\| _{k,l,\mu ,\alpha }^{\prime }$ by $\left\| \cdot
\right\| _{p,q,\mu ,\alpha }$ for some $p,q$ (and vice versa which is
completely analogous). But for all $f\in \mathcal{N}$ we have $\left|
f\right| _{k}^{\prime }\leq C\left| f\right| _{p}$ for some constant $C$ and
some $p$, since $\left| \cdot \right| _{k}^{\prime }$ has to be continuous
with respect to the projective limit topology on $\mathcal{N}$. That means
that the injection $i$ from $\mathcal{H}_{p}$ into the completion $\mathcal{K%
}_{k}$ of $\mathcal{N}$ with respect to $\left| \cdot \right| _{k}^{\prime }$
is a mapping bounded by $C$. We denote by $i$ also its linear extension from 
$\mathcal{H}_{p,\QTR{mathbb}{C}}$ into $\mathcal{K}_{k,\QTR{mathbb}{C}}$. It
follows that $i^{\otimes n}$ is bounded by $C^{n}$ from $\mathcal{H}_{p,%
\QTR{mathbb}{C}}^{\otimes n}$ into $\mathcal{K}_{k,\QTR{mathbb}{C}}^{\otimes
n}$. Now we choose $q$ such that $2^{\frac{q-l}{2}}\geq C$. Then 
\begin{eqnarray*}
\left\| \cdot \right\| _{k,l,\mu ,\alpha }^{\prime } &=&\sum_{n=0}^{\infty
}\left( n!\right) ^{2}2^{nl}\left| \cdot \right| _{k}^{\prime 2} \\
&\leq &\sum_{n=0}^{\infty }\left( n!\right) ^{2}2^{nl}C^{2n}\left| \cdot
\right| _{p}^{2} \\
&\leq &\left\| \cdot \right\| _{p,q,\mu ,\alpha }
\end{eqnarray*}
which is exactly what we need.\hfill $\blacksquare $

\begin{lemma}
\label{1eq3}There exist $p,C,K>0$ such that for all $n\in 
\QTR{mathbb}{N}_{0}$%
\begin{equation}
\int \left| P_{n}^{\mu ,\alpha }\left( z\right) \right| _{-p}^{2}d\mu \left(
z\right) \leq 4\left( n!\right) ^{2}C^{n}K.  \label{1eq4}
\end{equation}
\end{lemma}

\noindent \textbf{Proof. }We can use the estimate (\ref{1eq25}) and Lemma 
\ref{1eq11} to conclude the result.

\hfill $\blacksquare $

\begin{theorem}
There exists $p^{\prime },q^{\prime }>0$ such that for all $p\geq p^{\prime }
$, $q\geq q^{\prime }$ the topological embedding $(\mathcal{H}_{p})_{q,\mu
,\alpha }^{1}\subset L^{2}(\mu )$ holds.
\end{theorem}

\noindent \textbf{Proof. }Elements of the space $(\mathcal{N})_{\mu ,\alpha
}^{1}$ are defined as series convergent in the given topology. Now we need
the convergence of the series in $L^{2}(\mu )$. Choose $q^{\prime }$ such
that $C>2^{q^{\prime }}$ ($C$ from estimate (\ref{1eq4})). Let us
take an arbitrary 
\[
\varphi =\sum_{n=0}^{\infty }\left\langle P_{n}^{\mu ,\alpha },\varphi
_{\alpha }^{\left( n\right) }\right\rangle \in \mathcal{P}\left( \mathcal{N}%
^{\prime }\right) . 
\]
For $p>p^{\prime }$ ($p^{\prime }$ from the Lemma \ref{1eq3}) and 
$q>q^{\prime }$ the following estimates hold 
\begin{eqnarray*}
\left\| \varphi \right\| _{L^{2}(\mu )} &\leq &\sum_{n=0}^{\infty }\left\|
\left\langle P_{n}^{\mu ,\alpha },\varphi _{\alpha }^{\left( n\right)
}\right\rangle \right\| _{L^{2}(\mu )} \\
&\leq &\sum_{n=0}^{\infty }\left| \varphi _{\alpha }^{\left( n\right)
}\right| _{p}\left\| \left| P_{n}^{\mu ,\alpha }\right| _{-p}\right\|
_{L^{2}(\mu )} \\
&\leq &2K^{1/2}\sum_{n=0}^{\infty }n!2^{nq/2}\left| \varphi _{\alpha
}^{\left( n\right) }\right| _{p}\left( C2^{-q}\right) ^{n/2} \\
&\leq &2K^{1/2}\left( \sum_{n=0}^{\infty }\left( C2^{-q}\right) ^{n}\right)
^{1/2}\left( \sum_{n=0}^{\infty }\left( n!\right) ^{2}2^{nq}\left| \varphi
_{\alpha }^{\left( n\right) }\right| _{p}^{2}\right) ^{1/2} \\
&=&2K^{1/2}\left( 1-C2^{-q}\right) ^{-1/2}\left\| \varphi \right\| _{p,q,\mu
,\alpha }.
\end{eqnarray*}
Taking the closure the inequality extends to the whole space $(\mathcal{H}%
_{p})_{q,\mu ,\alpha }^{1}\,$.\hfill $\blacksquare $

\begin{corollary}
$(\mathcal{N})_{\mu ,\alpha }^{1}$ is continuously and densely embedded in $%
L^{2}(\mu )$.
\end{corollary}

\subsection{Description of test functions}

\begin{proposition}
\label{1eq5}Any test function $\varphi $ in $(\mathcal{N})_{\mu ,\alpha
}^{1}$ has a uniquely defined extension to $\mathcal{N}_{\QTR{mathbb}{C}%
}^{\prime }$ as an element of $\mathcal{E}_{\min }^{1}(\mathcal{N}_{%
\QTR{mathbb}{C}}^{\prime })$.
\end{proposition}

\noindent \textbf{Proof. }Any element $\varphi $ in $(\mathcal{N})_{\mu
,\alpha }^{1}$ is defined as a series of the following type 
\[
\varphi =\sum_{n=0}^{\infty }\left\langle P_{n}^{\mu ,\alpha },\varphi
_{\alpha }^{\left( n\right) }\right\rangle ,\quad \varphi _{\alpha }^{\left(
n\right) }\in \mathcal{N}_{\QTR{mathbb}{C}}^{\widehat{\otimes }n},
\]
such that 
\[
\left\| \varphi \right\| _{p,q,\mu ,\alpha }^{2}=\sum_{n=0}^{\infty }\left(
n!\right) ^{2}2^{nq}\left| \varphi _{\alpha }^{\left( n\right) }\right|
_{p}^{2}<\infty 
\]
for each $p,q\in \QTR{mathbb}{N}${. So we need to show the convergence of
the series} 
\[
\sum_{n=0}^{\infty }\left\langle P_{n}^{\mu ,\alpha }\left( z\right)
,\varphi _{\alpha }^{\left( n\right) }\right\rangle ,\quad z\in \mathcal{H}%
_{-p,\QTR{mathbb}{C}}
\]
to an entire function in $z$. Let $\epsilon >0$ and $\sigma _{\epsilon }>0$
as in (P$_{\alpha }$6) of Proposition \ref{1eq19}. We use (\ref{1eq25})
and estimate as follows
\begin{eqnarray*}
&&\sum_{n=0}^{\infty }\left| \left\langle P_{n}^{\mu ,\alpha }\left(
z\right) ,\varphi _{\alpha }^{\left( n\right) }\right\rangle \right|  \\
&\leq &\sum_{n=0}^{\infty }\left| P_{n}^{\mu ,\alpha }\left( z\right)
\right| _{-p}\left| \varphi _{\alpha }^{\left( n\right) }\right| _{p} \\
\  &\leq &2\sum_{n=0}^{\infty }n!\left| \varphi _{\alpha }^{\left( n\right)
}\right| _{p}\sigma _{\epsilon }^{-n} \\
\  &\leq &2\exp \left( \epsilon \left| z\right| _{-p^{\prime }}\right)
\left( \sum_{n=0}^{\infty }\left( n!\right) ^{2}2^{nq}\left| \varphi
_{\alpha }^{\left( n\right) }\right| _{p}^{2}\right) ^{1/2}\left(
\sum_{n=0}^{\infty }2^{-nq}\sigma _{\epsilon }^{-2n}\right) ^{1/2} \\
\  &\leq &2\left\| \varphi \right\| _{p,q,\mu ,\alpha }\left( 1-2^{-q}\sigma
_{\epsilon }^{-2}\right) ^{-1/2}\exp \left( \epsilon \left| z\right|
_{-p^{\prime }}\right) ,
\end{eqnarray*}
if $2^{q}>\sigma _{\epsilon }^{-2}$ and $p^{\prime }$ is such that $\mathcal{%
H}_{p}\hookrightarrow \mathcal{H}_{p^{\prime }}$ is Hilbert-Schmidt. That
means the series 
\[
\sum_{n=0}^{\infty }\left\langle P_{n}^{\mu ,\alpha }\left( z\right)
,\varphi _{\alpha }^{\left( n\right) }\right\rangle 
\]
converges uniformly and absolutely in any neighborhood of zero of any space $%
\mathcal{H}_{-p,\QTR{mathbb}{C}}$. Since each term $\langle P_{n}^{\mu
,\alpha }\left( z\right) ,\varphi _{\alpha }^{\left( n\right) }\rangle $ is
entire in $z$ the uniform convergence implies that 
\[
z\longmapsto \sum_{n=0}^{\infty }\left\langle P_{n}^{\mu ,\alpha }\left(
z\right) ,\varphi _{\alpha }^{\left( n\right) }\right\rangle 
\]
is entire on each $\mathcal{H}_{-p,\QTR{mathbb}{C}}$ and hence on $\mathcal{N%
}_{\QTR{mathbb}{C}}^{\prime }$. This complete the proof. \hfill $%
\blacksquare $

The following corollary gives an explicit estimate on the growth of test
functions and is a consequence of the above Proposition.

\begin{corollary}
\label{1eq6}For all $p>p^{\prime }$ such that the embedding $%
\mathcal{H}_{p}\hookrightarrow \mathcal{H}_{p^{\prime }}$ is of the
Hilbert-Schmidt class and for all $\epsilon >0$ there exists $\sigma
_{\epsilon }$ ($\sigma _{\epsilon }$ from Proposition \ref{1eq19}),
such that for $p\in \QTR{mathbb}{N}$ we obtain the following bound 
\[
\left| \varphi \left( z\right) \right| \leq C\left\| \varphi \right\|
_{p,q,\mu ,\alpha }\exp \left( \epsilon \left| z\right| _{-p^{\prime
}}\right) ,\;\varphi \in \left( \mathcal{N}\right) _{\mu ,\alpha
}^{1},\;z\in \mathcal{H}_{-p,\QTR{mathbb}{C}}\,,
\]
where $2^{q}>\sigma _{\epsilon }^{-2}$ and 
\[
C=2\left( 1-2^{-q}\sigma _{\epsilon }^{-2}\right) ^{-1/2}.
\]
\end{corollary}

\begin{remark}
Proposition \ref{1eq5} states 
\[
\left( \mathcal{N}\right) _{\mu ,\alpha }^{1}\subseteq \mathcal{E}_{\min
}^{1}\left( \mathcal{N}^{\prime }\right) 
\]
as sets, where 
\[
\mathcal{E}_{\min }^{1}\left( \mathcal{N}^{\prime }\right) =\left\{ \varphi
|_{\mathcal{N}^{\prime }}\,|\,\varphi \in \mathcal{E}_{\min }^{1}\left( 
\mathcal{N}_{\QTR{mathbb}{C}}^{\prime }\right) \right\} .
\]
Now we are going to show that the converse also holds.
\end{remark}

\begin{theorem}
\label{1eq7}For all functions $\alpha \in \mathrm{Hol}_{0}(\mathcal{N}_{%
\QTR{mathbb}{C}},\mathcal{N}_{\QTR{mathbb}{C}})$, as in Subsection \ref
{1eq15}, and for all measure $\mu \in \mathcal{M}_{a}(\mathcal{N}%
^{\prime })$, we have the topological identity 
\[
\left( \mathcal{N}\right) _{\mu ,\alpha }^{1}=\mathcal{E}_{\min }^{1}\left( 
\mathcal{N}^{\prime }\right) .
\]
\end{theorem}

\noindent \textbf{Proof. }Let $\varphi (z)\in \mathcal{E}_{\min }^{1}(%
\mathcal{N}^{\prime })$ be given such that 
\[
\varphi (z)=\sum_{n=0}^{\infty }\left\langle z^{\otimes n},\psi ^{\left(
n\right) }\right\rangle , 
\]
with 
\[
\left| \!\left| \!\left| \varphi \right| \!\right| \!\right|
_{p,q,1}^{2}=\sum_{n=0}^{\infty }\left( n!\right) ^{2}2^{nq}\left| \psi
^{\left( n\right) }\right| _{p}^{2}<\infty 
\]
for each $p,q\in \QTR{mathbb}{N}$. So we have 
\[
\left| \psi ^{\left( n\right) }\right| _{p}\leq \left( n!\right)
^{-1}2^{-nq/2}\left| \!\left| \!\left| \varphi \right| \!\right| \!\right|
_{p,q,1}. 
\]
On the other hand, we can use (\ref{1eq21}) to evaluate $\varphi (z)$ as 
\begin{eqnarray*}
\varphi (z) &=&\sum_{n=0}^{\infty }\left\langle z^{\otimes n},\psi ^{\left(
n\right) }\right\rangle \\
\ &=&\sum_{n=0}^{\infty }\left\langle \sum_{k=0}^{n}\sum_{m=0}^{k}\binom{n}{k%
}\frac{1}{m!}\left\langle P_{m}^{\mu ,\alpha }\left( z\right)
,B_{k}^{m}\right\rangle \widehat{\otimes }M_{n-k}^{\mu },\psi ^{\left(
n\right) }\right\rangle \\
\ &=&\sum_{n=0}^{\infty }\sum_{k=0}^{n}\sum_{m=0}^{k}\binom{n}{k}\frac{1}{m!}%
\left\langle \left\langle P_{m}^{\mu ,\alpha }\left( z\right)
,B_{k}^{m}\right\rangle ,\left( M_{n-k}^{\mu },\psi ^{\left( n\right)
}\right) _{\mathcal{H}^{\widehat{\otimes }\left( n-k\right) }}\right\rangle
\\
\ &=&\sum_{n=0}^{\infty }\sum_{k=0}^{n}\sum_{m=0}^{k}\binom{n}{k}\frac{1}{m!}
\left\langle P_{m}^{\mu ,\alpha }\left( z\right) ,\left\langle
B_{k}^{m},\left( M_{n-k}^{\mu },\psi ^{\left( n\right) }\right) _{\mathcal{H}%
^{\widehat{\otimes }\left( n-k\right) }}\right\rangle \right\rangle \\
\ &=&\sum_{m=0}^{\infty }\sum_{n=0}^{\infty }\sum_{k=0}^{n}\binom{n+m}{k+m}%
\frac{1}{m!}\left\langle P_{m}^{\mu ,\alpha }\left( z\right) ,\left\langle
B_{k+m}^{m},\left( M_{n-k}^{\mu },\psi ^{\left( n+m\right) }\right) _{%
\mathcal{H}^{\widehat{\otimes }\left( n-k\right) }}\right\rangle
\right\rangle \\
\ &=&\sum_{m=0}^{\infty }\left\langle P_{m}^{\mu ,\alpha }\left( z\right)
,\sum_{n=0}^{\infty }\sum_{k=0}^{n}\binom{n+m}{k+m}\frac{1}{m!}\left\langle
B_{k+m}^{m},\left( M_{n-k}^{\mu },\psi ^{\left( n+m\right) }\right) _{%
\mathcal{H}^{\widehat{\otimes }\left( n-k\right) }}\right\rangle
\right\rangle ,
\end{eqnarray*}
such that, if 
\[
\varphi (z)=\sum_{m=0}^{\infty }\left\langle P_{m}^{\mu ,\alpha }\left(
z\right) ,\varphi _{\alpha }^{\left( m\right) }\right\rangle , 
\]
then we conclude that 
\[
\varphi _{\alpha }^{\left( m\right) }=\sum_{n=0}^{\infty }\sum_{k=0}^{n}%
\binom{n+m}{k+m}\frac{1}{m!}\left\langle B_{k+m}^{m},\left( M_{n-k}^{\mu
},\psi ^{\left( n+m\right) }\right) _{\mathcal{H}^{\widehat{\otimes }\left(
n-k\right) }}\right\rangle . 
\]
Now for $p\in \QTR{mathbb}{N}$ we need estimate $|\varphi _{\alpha }^{\left(
n\right) }|_{p}$ by $\left| \!\left| \!\left| \cdot \right| \!\right|
\!\right| _{p,q,1}$ since the nuclear topology given by the norms $\left|
\!\left| \!\left| \cdot \right| \!\right| \!\right| _{p,q,1}$, is equivalent
to the projective topology induced by the norms $\mathrm{n}_{p,l,k}$ (see 
\cite{KSWY95}). Now we estimate $\varphi _{\alpha }^{\left( m\right) }$ as
follows 
\begin{eqnarray*}
\left| \varphi _{\alpha }^{\left( m\right) }\right| _{p} &\leq
&\sum_{n=0}^{\infty }\sum_{k=0}^{n}\binom{n+m}{k+m}\frac{1}{m!}\left|
B_{k+m}^{m}\right| _{\mathcal{H}_{-p}^{\widehat{\otimes }\left( k+m\right)
}\otimes \mathcal{H}_{p}^{\widehat{\otimes }m}}\left| \left( M_{n-k}^{\mu
},\psi ^{\left( n+m\right) }\right) _{\mathcal{H}^{\widehat{\otimes }\left(
n-k\right) }}\right| _{p} \\
\ &\leq &\sum_{n=0}^{\infty }\sum_{k=0}^{n}\binom{n+m}{k+m}\frac{1}{m!}%
\left| B_{k+m}^{m}\right| _{\mathcal{H}_{-p}^{\widehat{\otimes }\left(
k+m\right) }\otimes \mathcal{H}_{p}^{\widehat{\otimes }m}}\left|
M_{n-k}^{\mu }\right| _{-p}\left| \psi ^{\left( n+m\right) }\right| _{p}\,.
\end{eqnarray*}
Let us, at first, estimate the norm 
\[
\left| B_{k+m}^{m}\right| _{-p,p}:=\left| B_{k+m}^{m}\right| _{\mathcal{H}%
_{-p}^{\widehat{\otimes }\left( k+m\right) }\otimes \mathcal{H}_{p}^{%
\widehat{\otimes }m}}. 
\]
To do this we choose $p>p_{\mu }$ such that $\left\| i_{p,p_{\mu }}\right\|
_{HS}$ is finite and define 
\[
D_{\alpha ,\epsilon }:=\sup\limits_{\left| \theta \right| _{p}=\epsilon
}\left| g_{\alpha }\left( \theta \right) \right| _{p}\quad \mathrm{and\quad }%
\tilde{\epsilon}:=\frac{\epsilon }{e\left\| i_{p,p_{\mu }}\right\| _{HS}}. 
\]
So, with this 
\begin{eqnarray*}
\left| B_{m}^{n}\right| _{-p,p} &\leq &\sum_{l_{1},{\ldots },l_{n}=m}\frac{m!%
}{l_{1}!{\cdots }l_{n}!}\left| g_{\alpha }^{(l_{1})}\left( 0\right) \right|
_{-p,p}{\cdots }\left| g_{\alpha }^{(l_{n})}\left( 0\right) \right| _{-p,p}
\\
\ &\leq &\sum_{l_{1},{\ldots },l_{n}=m}\frac{m!l_{1}!{\cdots }l_{n}!}{l_{1}!{%
\cdots }l_{n}!}D_{\alpha ,\epsilon }^{n}\tilde{\epsilon}^{-m} \\
\ &\leq &m!D_{\alpha ,\epsilon }^{n}2^{m}\tilde{\epsilon}^{-m},
\end{eqnarray*}
that means 
\[
\left| B_{k+m}^{m}\right| _{-p,p}\leq \left( k+m\right) !D_{\alpha ,\epsilon
}^{m}2^{k+m}\tilde{\epsilon}^{-(k+m)}. 
\]
Now let $q\in \QTR{mathbb}{N}$ such that $2^{q/2}>K_{p}$ ($K_{p}:=eC\left\|
i_{p,p_{\mu }}\right\| _{HS}$ as in (\ref{1eq13})) and such that $2/(%
\tilde{\epsilon}K_{p})<1$, then we obtain 
\begin{eqnarray*}
&&\left| \varphi _{\alpha }^{\left( m\right) }\right| _{p} \\
&\leq &\sum_{n=0}^{\infty }\sum_{k=0}^{n}\binom{n+m}{k+m}\frac{1}{m!}\left(
m+k\right) !D_{\alpha ,\epsilon }^{m}\frac{2^{k+m}}{\tilde{\epsilon}^{k+m}}%
\left( n-k\right) !\left( K_{p}\right) ^{n-k}\frac{2^{-\left( n+m\right) q/2}%
}{\left( n+m\right) !}\left| \!\left| \!\left| \varphi \right| \!\right|
\!\right| _{p,q,1} \\
\ &\leq &\left| \!\left| \!\left| \varphi \right| \!\right| \!\right|
_{p,q,1}\frac{2^{-mq/2}}{m!}D_{\alpha ,\epsilon }^{m}\sum_{n=0}^{\infty
}\left( 2^{-q/2}K_{p}\right) ^{n}\sum_{k=0}^{n}\left( \frac{2}{\tilde{%
\epsilon}K_{p}}\right) ^{k} \\
\ &\leq &\left| \!\left| \!\left| \varphi \right| \!\right| \!\right|
_{p,q,1}\frac{2^{-mq/2}2^{m}}{m!\tilde{\epsilon}^{m}}D_{\alpha ,\epsilon
}^{m}\left( 1-2^{-q/2}K_{p}\right) ^{-1}\frac{\tilde{\epsilon}K_{p}}{\tilde{%
\epsilon}K_{p}-2} \\
\ &\equiv &L_{p,q,\alpha ,\tilde{\epsilon}}\frac{2^{-mq/2}2^{m}}{m!\tilde{%
\epsilon}^{m}}D_{\alpha ,\epsilon }^{m}\left| \!\left| \!\left| \varphi
\right| \!\right| \!\right| _{p,q,1}.
\end{eqnarray*}

\noindent For $q^{\prime }<q$ such that $2^{2}\tilde{\epsilon}%
^{-2}2^{(q^{\prime }-q)}D_{\alpha ,\epsilon }<1$ this follows the following
estimate 
\begin{eqnarray*}
\left\| \varphi \right\| _{p,q^{\prime },\mu ,\alpha }^{2} &\leq
&\sum_{m=0}^{\infty }\left( m!\right) ^{2}2^{mq^{\prime }}\left| \varphi
^{\left( m\right) }\right| _{p}^{2} \\
\ &\leq &\left| \!\left| \!\left| \varphi \right| \!\right| \!\right|
_{p,q,1}^{2}L_{p,q,\alpha ,\tilde{\epsilon}}^{2}\sum_{m=0}^{\infty }\left(
2^{2}\tilde{\epsilon}^{-2}2^{(q^{\prime }-q)}D_{\alpha ,\epsilon }\right)
^{m}<\infty .
\end{eqnarray*}
This complete the proof.\hfill $\blacksquare $

\bigskip \ Since we now have proved that the space of test functions $(%
\mathcal{N})_{\mu ,\alpha }^{1}$ is isomorphic to $\mathcal{E}_{\min }^{1}(%
\mathcal{N}^{\prime })$, for all measures $\mu \in \mathcal{M}_{a}(\mathcal{N%
}^{\prime })$ and for all holomorphic invertible function $\alpha \in 
\mathrm{Hol}_{0}(\mathcal{N}_{\QTR{mathbb}{C}},\mathcal{N}_{\QTR{mathbb}{C}%
}) $, such that $\alpha \left( 0\right) =0$, we will now drop the subscript $%
\mu ,\alpha $. The test function space $(\mathcal{N})^{1}$ is the same for
all measures and functions $\alpha $ in the above conditions.

\begin{corollary}
$(\mathcal{N})^{1}$ is an algebra under pointwise multiplication.
\end{corollary}

\begin{corollary}
$(\mathcal{N})^{1}$ admits `scaling', i.e., for $\lambda \in \QTR{mathbb}{C}$
the scaling operator $\sigma _{\lambda }:(\mathcal{N})^{1}\rightarrow (%
\mathcal{N})^{1}$ defined by $\sigma _{\lambda }\varphi \left( x\right)
:=\varphi \left( \lambda x\right) $, $\varphi \in (\mathcal{N})^{1}$, $x\in 
\mathcal{N}^{\prime }$ is well-defined.
\end{corollary}

\begin{corollary}
For all $z\in \mathcal{N}_{\QTR{mathbb}{C}}^{\prime }$ the space $(\mathcal{N%
})^{1}$ is invariant under the shift operator $\tau _{z}:\varphi \mapsto
\varphi \left( \cdot +z\right) $.
\end{corollary}

\section{Distributions}

In this section we will introduce and study the space $(\mathcal{N})_{\mu
,\alpha }^{-1}$ of distributions corresponding to the space of test
functions $(\mathcal{N})^{1}$ ($\equiv (\mathcal{N})_{\mu ,\alpha }^{1}$).
The goal is to prove that, for a fixed measure $\mu $ and for all function $%
\alpha $, as in the subsection \ref{1eq15}, the space $(\mathcal{N}%
)_{\mu ,\alpha }^{-1}=(\mathcal{N})_{\mu }^{-1}$, see Theorem \ref
{1eq41} below.

Since $\mathcal{P}(\mathcal{N}^{\prime })\subset (\mathcal{N})^{1}$ the
space $(\mathcal{N})_{\mu ,\alpha }^{-1}$ can be viewed as a subspace of $%
\mathcal{P}_{\mu }^{\prime }(\mathcal{N}^{\prime })$, i.e., 
\[
\left( \mathcal{N}\right) _{\mu ,\alpha }^{-1}\subset \mathcal{P}_{\mu
}^{\prime }\left( \mathcal{N}^{\prime }\right) . 
\]

\noindent Let us now introduce the Hilbert subspace $(\mathcal{H}%
_{-p})_{-q,\mu ,\alpha }^{-1}$ of $\mathcal{P}_{\mu }^{\prime }(\mathcal{N}%
^{\prime })$ for which the norm 
\[
\left\| \Phi \right\| _{-p,-q,\mu ,\alpha }^{2}:=\sum_{n=0}^{\infty
}2^{-qn}\left| \Phi _{\alpha }^{\left( n\right) }\right| _{-p}^{2} 
\]
is finite. Here we used the canonical representation 
\[
\Phi =\sum_{n=0}^{\infty }Q_{n}^{\mu ,\alpha }\left( \Phi _{\alpha }^{\left(
n\right) }\right) \in \mathcal{P}_{\mu }^{\prime }\left( \mathcal{N}^{\prime
}\right) 
\]
from Theorem \ref{1eq39}. The space $(\mathcal{H}_{-p})_{-q,\mu
,\alpha }^{-1}$ is the dual space of $(\mathcal{H}_{p})_{q,\mu ,\alpha }^{1}$
with respect to $L^{2}(\mu )$ (because of the biorthogonality of \textrm{%
\textbf{P}}$^{\mu ,\alpha }$- and \textrm{\textbf{Q}}$^{\mu ,\alpha }$%
-systems). By general duality theory 
\[
\left( \mathcal{N}\right) _{\mu ,\alpha }^{-1}=\bigcup_{p,q\in \QTR{mathbb}{N%
}}\left( \mathcal{H}_{-p}\right) _{-q,\mu ,\alpha }^{-1} 
\]
is the dual space of $\left( \mathcal{N}\right) ^{1}$ with respect to $%
L^{2}(\mu )$. As noted in Section \ref{1eq8} there exists a natural topology
on co-nuclear spaces (which coincide with the inductive limit topology). We
will consider $(\mathcal{N})_{\mu ,\alpha }^{-1}$ as a topological vector
space with this topology. So we have the nuclear triple 
\[
\left( \mathcal{N}\right) ^{1}\subset L^{2}\left( \mu \right) \subset \left( 
\mathcal{N}\right) _{\mu ,\alpha }^{-1}. 
\]

\noindent The action of a distribution 
\[
\Phi =\sum_{n=0}^\infty Q_n^{\mu ,\alpha }(\Phi _\alpha ^{\left( n\right)
})\in (\mathcal{N})_{\mu ,\alpha }^{-1} 
\]
on a test function 
\[
\varphi =\sum_{n=0}^\infty \langle P_n^{\mu ,\alpha },\varphi _\alpha
^{\left( n\right) }\rangle \in \left( \mathcal{N}\right) ^1 
\]
is given by 
\[
\left\langle \!\left\langle \Phi ,\varphi \right\rangle \!\right\rangle _\mu
=\sum_{n=0}^\infty n!\left\langle \Phi _\alpha ^{\left( n\right) },\varphi
_\alpha ^{\left( n\right) }\right\rangle . 
\]

For a more detailed characterization of the singularity of distributions in $%
(\mathcal{N})_{\mu ,\alpha }^{-1}$ we will introduce some subspaces in this
distribution space. For $\beta \in \left[ 0,1\right] $ we define 
\[
\begin{array}{ll}
\left( \mathcal{H}_{-p}\right) _{-q,\mu ,\alpha }^{-\beta } & :=\displaystyle%
\left\{ \Phi \in \mathcal{P}_{\mu }^{\prime }\left( \mathcal{N}^{\prime
}\right) \,|\,\sum_{n=0}^{\infty }\left( n!\right) ^{1-\beta }2^{-nq}\left|
\Phi _{\alpha }^{\left( n\right) }\right| _{-p}^{2}<\infty \right. \\ 
& \displaystyle\quad \quad \quad \quad \quad \quad \quad \quad \quad \quad
\quad \left. \mathrm{for}\;\Phi =\sum_{n=0}^{\infty }Q_{n}^{\mu ,\alpha
}\left( \Phi _{\alpha }^{\left( n\right) }\right) \right\}
\end{array}
\]
and 
\[
\left( \mathcal{N}\right) _{\mu ,\alpha }^{-\beta }=\bigcup_{p,q\in 
\QTR{mathbb}{N}}\left( \mathcal{H}_{-p}\right) _{-q,\mu ,\alpha }^{-\beta }. 
\]
It is clear that the singularity increases with increasing $\beta :$%
\[
\left( \mathcal{N}\right) _{\mu ,\alpha }^{-0}\subset \left( \mathcal{N}%
\right) _{\mu ,\alpha }^{-\beta _{1}}\subset \left( \mathcal{N}\right) _{\mu
,\alpha }^{-\beta _{2}}\subset \left( \mathcal{N}\right) _{\mu ,\alpha
}^{-1} 
\]
if $\beta _{1}\leq \beta _{2}$. We will also consider $(\mathcal{N})_{\mu
,\alpha }^{-\beta }$ as equipped with the natural topo\-logy.

\begin{example}
\emph{\textbf{(Generalized Radon-Nikodym derivative)}} We want to define a
generalized function $\rho _{\mu }^{\alpha }\left( z,\cdot \right) \in (%
\mathcal{N})_{\mu ,\alpha }^{-1}$, $z\in \mathcal{N}_{\QTR{mathbb}{C}%
}^{\prime }$ with the following property 
\[
\left\langle \!\left\langle \rho _{\mu }^{\alpha }\left( z,\cdot \right)
,\varphi \right\rangle \!\right\rangle _{\mu }=\int_{\mathcal{N}^{\prime
}}\varphi \left( x-z\right) d\mu \left( x\right) ,\quad \varphi \in \left( 
\mathcal{N}\right) ^{1}.
\]
That means we have to establish the continuity of $\rho _{\mu }^{\alpha
}\left( z,\cdot \right) $. Let $z\in \mathcal{H}_{-p,\QTR{mathbb}{C}}$. If $%
p\geq p^{\prime }\,$ is sufficiently large and $\epsilon >0$ small enough,
Corollary \ref{1eq6} applies, i.e., $\exists q\in \QTR{mathbb}{N}$
and $C>0$ such that 
\begin{eqnarray*}
\left| \int_{\mathcal{N}^{\prime }}\varphi \left( x-z\right) d\mu \left(
x\right) \right|  &\leq &C\left\| \varphi \right\| _{p,q,\mu ,\alpha }\int_{%
\mathcal{N}^{\prime }}\exp \left( \epsilon \left| x-z\right| _{-p^{\prime
}}\right) d\mu \left( x\right)  \\
&\leq &C\left\| \varphi \right\| _{p,q,\mu ,\alpha }\exp \left( \epsilon
\left| z\right| _{-p^{\prime }}\right) \int_{\mathcal{N}^{\prime }}\exp
\left( \epsilon \left| x\right| _{-p^{\prime }}\right) d\mu \left( x\right) .
\end{eqnarray*}
If $\epsilon $ is chosen sufficiently small the last integral exists
(Lemma \ref{1eq11}-3). Thus we have in fact $\rho _{\mu }^{\alpha
}\left( z,\cdot \right) \in (\mathcal{N})_{\mu ,\alpha }^{-1}$. It is clear
that whenever the Radon-Nikodym derivative $\frac{d\mu (x+\xi )}{d\mu \left(
x\right) }$ exists (e.g., $\xi \in \mathcal{N}$ in case $\mu $ is $\mathcal{N%
}$-quasi-invariant) it coincides with $\rho _{\mu }^{\alpha }(\xi ,\cdot )$
defined above. We will show that in $(\mathcal{N})_{\mu ,\alpha }^{-1}$ we
have the canonical expansion 
\[
\rho _{\mu }^{\alpha }\left( z,\cdot \right) =\sum_{n=0}^{\infty }\frac{1}{n!%
}\left( -1\right) ^{n}\left\langle Q_{n}^{\mu ,\alpha }\left( \cdot \right)
,P_{n}^{\delta _{0},\alpha }\left( -z\right) \right\rangle 
\]
where $P_{n}^{\delta _{0},\alpha }\left( -z\right) $ is defined in (\ref
{1eq31}). It is easy to see that the r.h.s. defines an element in $(%
\mathcal{N})_{\mu ,\alpha }^{-1}$. Since both sides are in $(\mathcal{N}%
)_{\mu ,\alpha }^{-1}$ it is sufficient to compare their action on a total
set from $(\mathcal{N})^{1}$. For $\varphi _{\alpha }^{(n)}\in \mathcal{N}_{%
\QTR{mathbb}{C}}^{\widehat{\otimes }n}$ we have 
\begin{eqnarray*}
&&\left\langle \!\left\langle \rho _{\mu }^{\alpha }\left( z,\cdot \right)
,\left\langle P_{n}^{\mu ,\alpha }\left( \cdot \right) ,\varphi _{\alpha
}^{(n)}\right\rangle \right\rangle \!\right\rangle _{\mu } \\
&=&\left\langle \!\!\left\langle \sum_{k=0}^{\infty }\frac{1}{k!}\left(
-1\right) ^{k}\left\langle Q_{k}^{\mu ,\alpha }\left( \cdot \right)
,P_{k}^{\delta _{0},\alpha }\left( -z\right) \right\rangle ,\left\langle
P_{n}^{\mu ,\alpha }\left( \cdot \right) ,\varphi _{\alpha
}^{(n)}\right\rangle \right\rangle \!\!\right\rangle _{\mu } \\
&=&\left\langle P_{n}^{\delta _{0},\alpha }\left( -z\right) ,\varphi
_{\alpha }^{(n)}\right\rangle ,
\end{eqnarray*}
where we have used the biorthogonality property of the \textrm{\textbf{Q}}$%
^{\mu ,\alpha }$- and \textrm{\textbf{P}}$^{\mu ,\alpha }$- systems. On the
other hand 
\begin{eqnarray*}
&&\left\langle \!\left\langle \rho _{\mu }^{\alpha }\left( z,\cdot \right)
,\left\langle P_{n}^{\mu ,\alpha }\left( \cdot \right) ,\varphi _{\alpha
}^{(n)}\right\rangle \right\rangle \!\right\rangle _{\mu } \\
&=&\int_{\mathcal{N}^{\prime }}\left\langle P_{n}^{\mu ,\alpha }\left(
x-z\right) ,\varphi _{\alpha }^{(n)}\right\rangle d\mu \left( x\right)  \\
&=&\sum_{k=0}^{n}\binom{n}{k}\int_{\mathcal{N}^{\prime }}\left\langle
P_{k}^{\mu ,\alpha }\left( x\right) \widehat{\otimes }P_{n-k}^{\delta
_{0},,\alpha }\left( -z\right) ,\varphi _{\alpha }^{(n)}\right\rangle d\mu
\left( x\right)  \\
&=&\sum_{k=0}^{n}\binom{n}{k}\QTR{mathbb}{E}_{\mu }\left( \left\langle
P_{k}^{\mu ,\alpha }\left( \cdot \right) \widehat{\otimes }P_{n-k}^{\delta
_{0},,\alpha }\left( -z\right) ,\varphi _{\alpha }^{(n)}\right\rangle
\right)  \\
&=&\left\langle P_{n}^{\delta _{0},\alpha }\left( -z\right) ,\varphi
_{\alpha }^{(n)}\right\rangle ,
\end{eqnarray*}
where we made use of the relation (\ref{1eq24}). This had to be shown. In
other words, we have proven that $\rho _{\mu }^{\alpha }\left( z,\cdot
\right) $ is the generating function of the \textrm{\textbf{Q}}$^{\mu
,\alpha }$-system. 
\[
\rho _{\mu }^{\alpha }\left( -z,\cdot \right) =\sum_{n=0}^{\infty }\frac{1}{%
n!}\left\langle Q_{n}^{\mu ,\alpha }\left( \cdot \right) ,P_{n}^{\delta
_{0},\alpha }\left( z\right) \right\rangle 
\]
\end{example}

\begin{example}
\emph{\textbf{(Delta function)}} For $z\in \mathcal{N}_{\QTR{mathbb}{C}%
}^{\prime }$ we define a distribution by the following \textrm{\textbf{Q}}$%
^{\mu ,\alpha }$-decomposition: 
\[
\delta _{z}=\sum_{n=0}^{\infty }\frac{1}{n!}Q_{n}^{\mu ,\alpha }\left(
P_{n}^{\mu ,\alpha }\left( z\right) \right) .
\]
If $p\in \QTR{mathbb}{N}$ is large enough and $\epsilon >0$ sufficiently
small there exists $\sigma _{\epsilon }>0$ according to (\ref{1eq25}) such
that 
\begin{eqnarray*}
\left\| \delta _{z}\right\| _{-p,-q,\mu ,\alpha }^{2}=\sum_{n=0}^{\infty
}\left( n!\right) ^{-2}2^{-nq}\left| P_{n}^{\mu ,\alpha }\left( z\right)
\right| _{-p}^{2} \\
\leq 4\exp \left( 2\epsilon \left| z\right| _{-p}\right) \sum_{n=0}^{\infty
}\sigma _{\epsilon }^{-2n}2^{-nq},\quad \quad z\in \mathcal{H}_{-p,%
\QTR{mathbb}{C}},
\end{eqnarray*}
which is finite for sufficiently large $q\in \QTR{mathbb}{N}${.} Thus $%
\delta _{z}\in (\mathcal{N})_{\mu ,\alpha }^{-1}$.

For 
\[
\varphi =\sum_{n=0}^{\infty }\left\langle P_{n}^{\mu ,\alpha },\varphi
_{\alpha }^{\left( n\right) }\right\rangle \in (\mathcal{N})^{1}
\]
the action of $\delta _{z}$ is given by 
\[
\left\langle \!\left\langle \delta _{z},\varphi \right\rangle
\!\right\rangle _{\mu }=\sum_{n=0}^{\infty }\left\langle P_{n}^{\mu ,\alpha
}\left( z\right) ,\varphi _{\alpha }^{\left( n\right) }\right\rangle
=\varphi \left( z\right) 
\]
because of the biorthogonality property, see Theorem \ref{1eq37} pag.~%
\pageref{1eq37}. This means that $\delta _{z}$ (in particular for $z$
real) plays the role of a ``$\delta $-function'' (evaluation map) in the
calculus we discuss.
\end{example}

\begin{theorem}
\label{1eq41}For a fixed measure $\mu $ and for all function $%
\alpha $, as in subsection \ref{1eq15}, we have 
\[
\left( \mathcal{N}\right) _{\mu ,\alpha }^{-1}=\left( \mathcal{N}\right)
_{\mu }^{-1},
\]
i.e., the space of distributions is the same for all functions $\alpha $ in
the above conditions.
\end{theorem}

\noindent \textbf{Proof. }Let $\Phi \in (\mathcal{N})_{\mu ,\alpha }^{-1}$
be given, then by Theorem \ref{1eq39} there exists gene\-ralized
kernels $\Phi _{\alpha }^{\left( n\right) }\in \mathcal{N}_{\QTR{mathbb}{C}%
}^{\prime \widehat{\otimes }n}$ such that $\Phi $ has the following
representation 
\[
\Phi =\sum_{n=0}^{\infty }\left\langle Q_{n}^{\mu ,\alpha },\Phi _{\alpha
}^{\left( n\right) }\right\rangle .
\]
Now we use the definition of $Q_{n}^{\mu ,\alpha }$ given in (\ref
{1eq34}) to obtain 
\begin{eqnarray}
S_{\mu }\Phi \left( \theta \right)  &=&\sum_{n=0}^{\infty }\left\langle \Phi
_{\alpha }^{\left( n\right) },g_{\alpha }\left( \theta \right) ^{\otimes
n}\right\rangle   \nonumber \\
&=&S_{\mu }\widehat{\Phi }\left( g_{\alpha }\left( \theta \right) \right)
,\qquad \theta \in \mathcal{N}_{\QTR{mathbb}{C}},  \label{1eq42}
\end{eqnarray}
where 
\[
\widehat{\Phi }=\sum_{n=0}^{\infty }\left\langle Q_{n}^{\mu },\Phi _{\alpha
}^{\left( n\right) }\right\rangle \in (\mathcal{N})_{\mu }^{-1}.
\]
Hence by characterization Theorem \ref{1eq2} $S_{\mu }\widehat{%
\Phi }\in \mathrm{Hol}_{0}(\mathcal{N}_{\QTR{mathbb}{C}})$. But from (\ref
{1eq42}) we see that 
\[
S_{\mu }\Phi =\left( S_{\mu }\widehat{\Phi }\right) \circ g_{\alpha }\in 
\mathrm{Hol}_{0}\left( \mathcal{N}_{\QTR{mathbb}{C}}\right) ,
\]
since this is the composition of two holomorphic functions (see \cite{Di81}%
), again by the characterization Theorem \ref{1eq2} we
conclude that $\Phi \in (\mathcal{N})_{\mu }^{-1}$. Hence $(\mathcal{N}%
)_{\mu ,\alpha }^{-1}\subseteq (\mathcal{N})_{\mu }^{-1}$.

Conversely, let $\Psi \in (\mathcal{N})_{\mu }^{-1}$ be given, i.e., 
\[
\Psi =\sum_{n=0}^{\infty }\left\langle Q_{n}^{\mu },\Psi ^{\left( n\right)
}\right\rangle ,\quad \Psi ^{\left( n\right) }\in \mathcal{N}_{\QTR{mathbb}{C%
}}^{\prime \widehat{\otimes }n}.
\]
We want to prove that $\Psi \in (\mathcal{N})_{\mu ,\alpha }^{-1}$. Due to (%
\ref{1eq34}) and the definition of $(\mathcal{N})_{\mu }^{-1}$ it
is sufficient to show that 
\[
S_{\mu }\Psi \left( \theta \right) =\sum_{n=0}^{\infty }\left\langle 
\widehat{\Psi }_{\alpha }^{\left( n\right) },g_{\alpha }\left( \theta
\right) ^{\otimes n}\right\rangle ,\quad \theta \in \mathcal{N}_{%
\QTR{mathbb}{C}},
\]
where $\widehat{\Psi }_{\alpha }^{\left( n\right) }$ satisfy, for $p,q\in 
\QTR{mathbb}{N}$%
\[
\sum_{n=0}^{\infty }2^{-nq}\left| \widehat{\Psi }_{\alpha }^{\left( n\right)
}\right| _{-p}^{2}<\infty .
\]

\noindent On the other hand, for a given $\theta \in \mathcal{N}_{%
\QTR{mathbb}{C}}$%
\[
S_{\mu }\Psi \left( \theta \right) =\sum_{n=0}^{\infty }\left\langle \Psi
^{\left( n\right) },\theta ^{\otimes n}\right\rangle =:G\left( \theta
\right) 
\]
and, consequently $G\in \mathrm{Hol}_{0}(\mathcal{N}_{\QTR{mathbb}{C}})$.
But we can write 
\[
G\left( \theta \right) =G\left( \alpha \circ g_{\alpha }\left( \theta
\right) \right) =\widehat{G}(g_{\alpha }\left( \theta \right) ),
\]
where $\widehat{G}=G\circ \alpha $\thinspace , with $G\circ \alpha \in 
\mathrm{Hol}_{0}(\mathcal{N}_{\QTR{mathbb}{C}})$. Therefore 
\[
\widehat{G}(g_{\alpha }\left( \theta \right) )=\sum_{n=0}^{\infty
}\left\langle \widehat{G}_{\alpha }^{\left( n\right) },g_{\alpha }\left(
\theta \right) ^{\otimes n}\right\rangle ,
\]
where the coefficients $\widehat{G}_{\alpha }^{\left( n\right) }$ verify 
\[
\sum_{n=0}^{\infty }2^{-nq}\left| \widehat{G}_{\alpha }^{\left( n\right)
}\right| _{-p}^{2}<\infty .
\]
Therefore with $\widehat{\Psi }_{\alpha }^{\left( n\right) }=\widehat{G}%
_{\alpha }^{\left( n\right) }$ follows the result, i.e., $\Psi \in (\mathcal{%
N})_{\mu ,\alpha }^{-1}$.\hfill $\blacksquare $

\section{The Wick product}

Here we give the natural generalization of the \textbf{Wick multiplication}
in the present setting.

\begin{definition}
Let $\Phi ,\Psi \in (\mathcal{N})_{\mu }^{-1}.$ Then we define the \emph{%
\textbf{Wick product}} $\Phi \diamond \Psi $ by 
\[
S_{\mu }\left( \Phi \diamond \Psi \right) =S_{\mu }\Phi \cdot S_{\mu }\Psi
\,.
\]
\end{definition}

This is well defined because Hol$_{0}(\mathcal{N}_{\QTR{mathbb}{C}})$ is an
algebra and thus by characterization theorem there exists an element in $(%
\mathcal{N})_{\mu }^{-1}$ $\Phi \diamond \Psi $ such that $S_{\mu }(\Phi
\diamond \Psi )=S_{\mu }\Phi \cdot S_{\mu }\Psi .$

From this it follows 
\[
Q_{n}^{\mu ,\alpha }\left( \Phi _{\alpha }^{\left( n\right) }\right)
\diamond Q_{m}^{\mu ,\alpha }\left( \Psi _{\alpha }^{\left( m\right)
}\right) =Q_{n+m}^{\mu ,\alpha }\left( \Phi _{\alpha }^{\left( n\right) }%
\widehat{\otimes }\Psi _{\alpha }^{\left( m\right) }\right) , 
\]
$\Phi _{\alpha }^{\left( n\right) }\in \mathcal{N}_{\QTR{mathbb}{C}}^{\prime 
\widehat{\otimes }n}$ and $\Psi _{\alpha }^{\left( m\right) }\in \mathcal{N}%
_{\QTR{mathbb}{C}}^{\prime \widehat{\otimes }m}.$ So in terms of \textrm{%
\textbf{Q}}$^{\mu ,\alpha }$-decomposition 
\[
\Phi =\sum_{n=0}^{\infty }Q_{n}^{\mu ,\alpha }\left( \Phi _{\alpha }^{\left(
n\right) }\right) \quad \mathrm{and\quad }\Psi =\sum_{m=0}^{\infty
}Q_{m}^{\mu ,\alpha }\left( \Psi _{\alpha }^{\left( m\right) }\right) 
\]
the Wick product is given by 
\[
\Phi \diamond \Psi =\sum_{n=0}^{\infty }Q_{n}^{\mu ,\alpha }\left( \Xi
_{\alpha }^{\left( n\right) }\right) , 
\]
where 
\[
\Xi _{\alpha }^{\left( n\right) }=\sum_{k=0}^{n}\Phi _{\alpha }^{\left(
k\right) }\widehat{\otimes }\Psi _{\alpha }^{\left( n-k\right) }. 
\]
This allows for a concrete norm estimate.

\begin{proposition}
The Wick product is continuous on $(\mathcal{N})_{\mu }^{-1}.$ In particular
the following estimate holds for $\Phi \in (\mathcal{H}_{-p_{1}})_{-q_{1},%
\mu ,\alpha }^{-1}\,,$ $\Psi \in (\mathcal{H}_{-p_{2}})_{-q_{2},\mu ,\alpha
}^{-1}$ and $p=\max (p_{1},p_{2}),$ $q=q_{1}+q_{2}+1$%
\[
\left\| \Phi \diamond \Psi \right\| _{-p,-q,\mu ,\alpha }\leq \left\| \Phi
\right\| _{-p_{1},-q_{1},\mu ,\alpha }\left\| \Psi \right\|
_{-p_{2},-q_{2},\mu ,\alpha }\,. 
\]
\end{proposition}

\noindent \textbf{Proof. }We can estimate as follows 
\begin{eqnarray*}
\left\| \Phi \diamond \Psi \right\| _{-p,-q,\mu ,\alpha }^{2}
&=&\sum_{n=0}^{\infty }2^{-nq}\left| \Xi _{\alpha }^{\left( n\right)
}\right| _{-p}^{2} \\
&=&\sum_{n=0}^{\infty }2^{-nq}\left( \sum_{k=0}^{n}\left| \Phi _{\alpha
}^{\left( k\right) }\right| _{-p}\left| \Psi _{\alpha }^{\left( n-k\right)
}\right| _{-p}\right) ^{2} \\
&\leq &\sum_{n=0}^{\infty }2^{-nq}\left( n+1\right) \sum_{k=0}^{n}\left|
\Phi _{\alpha }^{\left( k\right) }\right| _{-p}^{2}\left| \Psi _{\alpha
}^{\left( n-k\right) }\right| _{-p}^{2} \\
&\leq &\sum_{n=0}^{\infty }\sum_{k=0}^{n}2^{-nq_{1}}\left| \Phi _{\alpha
}^{\left( k\right) }\right| _{-p}^{2}2^{-nq_{2}}\left| \Psi _{\alpha
}^{\left( n-k\right) }\right| _{-p}^{2} \\
&\leq &\left( \sum_{n=0}^{\infty }2^{-nq_{1}}\left| \Phi _{\alpha }^{\left(
n\right) }\right| _{-p}^{2}\right) \left( \sum_{n=0}^{\infty
}2^{-nq_{2}}\left| \Psi _{\alpha }^{\left( n\right) }\right| _{-p}^{2}\right)
\\
&=&\left\| \Phi \right\| _{-p_{1},-q_{1},\mu ,\alpha }^{2}\left\| \Psi
\right\| _{-p_{2},-q_{2},\mu ,\alpha }^{2}.
\end{eqnarray*}
\hfill $\blacksquare $

Similar to the Gaussian case the special properties of the space $(\mathcal{N%
})_\mu ^{-1}$ allow the definition of \textbf{Wick analytic functions}%
\textit{\ }under very general assum\-ptions. This has proven to be of some
relevance to solve equations e.g., of the type $\Phi \diamond X=\Psi $ for $%
X\in (\mathcal{N})_\mu ^{-1}$. See \cite{KLS96} for the Gaussian case.

\begin{proposition}
For any $n\in \QTR{mathbb}{N}$ and any $\alpha $ as in Subsection \ref
{1eq15} we have $Q_{n}^{\mu ,\alpha }=(Q_{1}^{\mu ,\alpha })^{\diamond
n}.$
\end{proposition}

\noindent \textbf{Proof. }Let $\Phi ^{\left( 1\right) }\in \mathcal{N}_{%
\QTR{mathbb}{C}}^{\prime }$ be given. Thus, if $\theta \in \mathcal{N}_{%
\QTR{mathbb}{C}},$ follows 
\begin{eqnarray*}
S_{\mu }\left[ \left( Q_{1}^{\mu ,\alpha }\left( \Phi ^{\left( 1\right)
}\right) \right) ^{\diamond n}\right] \left( \theta \right) &=&\left\langle
\Phi ^{\left( 1\right) },g_{\alpha }\left( \theta \right) \right\rangle ^{n}
\\
\ &=&\left\langle \left( \Phi ^{\left( 1\right) }\right) ^{\widehat{\otimes }%
n},\left( g_{\alpha }\left( \theta \right) \right) ^{\otimes n}\right\rangle
\\
\ &=&S_{\mu }\left[ Q_{n}^{\mu ,\alpha }\left( \left( \Phi ^{\left( 1\right)
}\right) ^{\widehat{\otimes }n}\right) \right] \left( \theta \right) .
\end{eqnarray*}
\hfill $\blacksquare $

\begin{theorem}
Let $F:\QTR{mathbb}{C}{\rightarrow \QTR{mathbb}{C}}$ be analytic in a
neighborhood of the point $z_{0}=\QTR{mathbb}{E}\left( \Phi \right) ,$ $\Phi
\in (\mathcal{N})_{\mu }^{-1}.$ Then $F^{\diamond }\left( \Phi \right) $
defined by $S_{\mu }\left( F^{\diamond }\left( \Phi \right) \right)
=F(S_{\mu }\Phi )$ exists in $(\mathcal{N})_{\mu }^{-1}.$
\end{theorem}

\noindent \textbf{Proof. }By Theorems \ref{1eq41} and \ref
{1eq2} we have $S_{\mu }\Phi \in \mathrm{Hol}_{0}(\mathcal{N}_{%
\QTR{mathbb}{C}}).$ Then $F(S_{\mu }\Phi )\in \mathrm{Hol}_{0}(\mathcal{N}_{%
\QTR{mathbb}{C}})$ since the composition of two analytic functions is also
analytic. Again by the above mentioned theorems we find that $F^{\diamond
}\left( \Phi \right) $ exists in $(\mathcal{N})_{\mu }^{-1}$. 
$\blacksquare $

\bigskip

\begin{remark}
If $F\left( z\right) $ have the following representation 
\[
F\left( z\right) =\sum_{n=0}^{\infty }a_{n}\left( z-z_{0}\right) ^{n},
\]
then the \textit{Wick series} 
\[
\sum_{n=0}^{\infty }a_{n}\left( \Phi -z_{0}\right) ^{\diamond n}
\]
(where $\Psi ^{\diamond n}=\Psi \diamond \cdots \diamond \Psi $ n-times)
converges in $(\mathcal{N})_{\mu }^{-1}$ and 
\[
F^{\diamond }\left( \Phi \right) =\sum_{n=0}^{\infty }a_{n}\left( \Phi
-z_{0}\right) ^{\diamond n}
\]
holds.
\end{remark}

\begin{example}
The above mentioned equation $\Phi \diamond X=\Psi $ can be solved if 
$\QTR{mathbb}{E}_{\mu }\left( \Phi \right) =S_{\mu }\Phi \left(
0\right) \neq 0.$ That implies $(S_{\mu }\Phi )^{-1}\in \mathrm{Hol}_{0}(%
\mathcal{N}_{\QTR{mathbb}{C}})$. Thus 
\[
\Phi ^{\diamond \left( -1\right) }=S_{\mu }^{-1}((S_{\mu }\Phi )^{-1})\in (%
\mathcal{N})_{\mu }^{-1}.
\]
Then $X=$ $\Phi ^{\diamond \left( -1\right) }\diamond \Psi \,$ is the
solution in $(\mathcal{N})_{\mu }^{-1}.$ For more instructive examples we
refer the reader to Section 5 of \cite{KLS96}.\vspace{1cm}
\end{example}

\noindent \textbf{Acknowledgements}

We are grateful to Dr.~N.~A.~Kachanovsky for helpful discussions. The first
author was supported by the DFG (Deutsche Forschungsgemeinschaft). J.L.S
acknowledges financial support of CITMA (Centro de Ci\^{e}ncias e Tecnologia
da Madeira).

\section{Change of measure \label{t75}}

Suppose we are given two measures $\mu ,\tilde{\mu}\in \mathcal{M}_{a}(%
\mathcal{N}^{\prime })$ both satisfying Assumption~\ref{t37}. Let a
distribution $\tilde{\Phi}\in (\mathcal{N})_{\tilde{\mu}}^{-1}$ be given.
Since the test function space $(\mathcal{N})^{1}$ is invariant under changes
of measure in view of Theorem~\ref{1eq7}, the continuous mapping 
\[
\varphi \longmapsto \langle \!\langle \tilde{\Phi},\varphi \rangle \!\rangle
_{\tilde{\mu}},\;\varphi \in (\mathcal{N})^{1},
\]
can also be represented as a distribution $\Phi \in (\mathcal{N})_{\mu }^{-1}
$. So we have the implicit relation 
\[
\tilde{\Phi}\in (\mathcal{N})_{\tilde{\mu}}^{-1}\longleftrightarrow \Phi \in
(\mathcal{N})_{\mu }^{-1},
\]
defined by 
\[
\langle \!\langle \tilde{\Phi},\varphi \rangle \!\rangle _{\tilde{\mu}%
}=\langle \!\langle \Phi ,\varphi \rangle \!\rangle _{\tilde{\mu}}.
\]
This section provide formulas which make this relation more explicit in
terms of re-decomposition of the $\QTR{mathbb}{Q}^{\mu ,\alpha }$-system.
First we need an explicit relation of the corresponding $\QTR{mathbb}{P}%
^{\mu ,\alpha }$-system.

\begin{lemma}
Let $\mu ,\tilde{\mu}\in \mathcal{M}_{a}(\mathcal{N}^{\prime })$ be given,
then 
\begin{equation}
P_{n}^{\mu ,\alpha }(x)=\sum_{k+m+l=n}\frac{n!}{k!m!l!}P_{k}^{\tilde{\mu}%
,\alpha }(x)\hat{\otimes}P_{m}^{\mu ,\alpha }(0)\hat{\otimes}M_{l}^{\tilde{%
\mu},\alpha }.  \label{eq1}
\end{equation}
\end{lemma}

\noindent \textbf{Proof. }Expanding each factor in the formula 
\[
e_{\mu }^{\alpha }\left( \theta ;x\right) =e_{\tilde{\mu}}^{\alpha }\left(
\theta ;x\right) l_{\mu }^{\alpha -1}\left( \theta \right) l_{\tilde{\mu}%
}^{\alpha }\left( \theta \right) , 
\]
we obtain 
\begin{eqnarray*}
&&\sum_{n=0}^{\infty }\frac{1}{n!}\langle P_{n}^{\mu ,\alpha }(x),\theta
^{\otimes n}\rangle \\
&=&\sum_{k=0}^{\infty }\frac{1}{k!}\langle P_{k}^{\tilde{\mu},\alpha
}(x),\theta ^{\otimes k}\rangle \sum_{m=0}^{\infty }\frac{1}{m!}\langle
P_{m}^{\mu ,\alpha }(0),\theta ^{\otimes m}\rangle \sum_{l=0}^{\infty }\frac{%
1}{l!}\langle M_{l}^{\tilde{\mu},\alpha },\theta ^{\otimes l}\rangle \\
&=&\sum_{n=0}^{\infty }\frac{1}{n!}\left\langle \sum_{k+m+l=n}\frac{n!}{%
k!m!l!}P_{k}^{\tilde{\mu},\alpha }\left( x\right) \hat{\otimes}P_{m}^{\mu
,\alpha }(0)\hat{\otimes}M_{l}^{\tilde{\mu},\alpha },\theta ^{\otimes
n}\right\rangle .
\end{eqnarray*}
A comparison of coefficients gives the above result.\hfill $\blacksquare $

An immediate consequence is the next reordering lemma.

\begin{lemma}
Let $\varphi \in (\mathcal{N})^{1}$ be given. Then $\varphi $ has the
representation in $\QTR{mathbb}{P}^{\mu ,\alpha }$-system as well as $%
\QTR{mathbb}{P}^{\tilde{\mu},\alpha }$-system: 
\[
\varphi =\sum_{n=0}^{\infty }\langle P_{n}^{\mu ,\alpha },\varphi _{\alpha
}^{(n)}\rangle =\sum_{n=0}^{\infty }\langle P_{n}^{\tilde{\mu},\alpha },%
\tilde{\varphi}_{\alpha }^{(n)}\rangle , 
\]
where $\varphi _{\alpha }^{(n)}$, $\tilde{\varphi}_{\alpha }^{(n)}\in 
\mathcal{N}_{\QTR{mathbb}{C}}^{\hat{\otimes}n}$ for all $n\in \QTR{mathbb}{N}%
_{0}$ and the following formula holds 
\begin{equation}
\tilde{\varphi}_{\alpha }^{(n)}=\sum_{m,l=0}^{\infty }\frac{(n+m+l)!}{n!m!l!}%
(P_{m}^{\mu ,\alpha }(0)\hat{\otimes}M_{l}^{\tilde{\mu},\alpha },\varphi
_{\alpha }^{(n+m+l)})_{\mathcal{H}^{\hat{\otimes}(m+l)}}.  \label{eq2}
\end{equation}
\end{lemma}

\noindent \textbf{Proof. }We use the relation (\ref{eq1}) to obtain 
\begin{eqnarray*}
\varphi &=&\sum_{n=0}^{\infty }\langle P_{n}^{\mu ,\alpha },\varphi _{\alpha
}^{(n)}\rangle \\
&=&\sum_{n=0}^{\infty }\left\langle \sum_{k+m+l=n}\frac{n!}{k!m!l!}P_{k}^{%
\tilde{\mu},\alpha }(x)\hat{\otimes}P_{m}^{\mu ,\alpha }(0)\hat{\otimes}%
M_{l}^{\tilde{\mu},\alpha },\varphi _{\alpha }^{(n)}\right\rangle \\
&=&\sum_{k,m,l=0}^{\infty }\frac{(k+m+l)!}{k!m!l!}\langle P_{k}^{\tilde{\mu}%
,\alpha }(x),(P_{m}^{\mu ,\alpha }(0)\hat{\otimes}M_{l}^{\tilde{\mu},\alpha
},\varphi _{\alpha }^{(k+m+l)})_{\mathcal{H}^{\hat{\otimes}(m+l)}}\rangle \\
&=&\sum_{k=0}^{\infty }\left\langle P_{k}^{\tilde{\mu},\alpha
}(x),\sum_{m,l=0}^{\infty }\frac{(k+m+l)!}{k!m!l!}(P_{m}^{\mu ,\alpha }(0)%
\hat{\otimes}M_{l}^{\tilde{\mu},\alpha },\varphi _{\alpha }^{(k+m+l)})_{%
\mathcal{H}^{\hat{\otimes}(m+l)}}\right\rangle .
\end{eqnarray*}
Then a comparison of coefficients give the result.\hfill $\blacksquare $

Now we may prove the announced theorem.

\begin{theorem}
Let $\tilde{\Phi}$ be a generalized function with representation 
\[
\tilde{\Phi}=\sum_{n=0}^{\infty }\langle Q_{n}^{\tilde{\mu},\alpha },\tilde{%
\Phi}_{\alpha }^{(n)}\rangle . 
\]
Then 
\[
\Phi =\sum_{n=0}^{\infty }\langle Q_{n}^{\mu ,\alpha },\Phi _{\alpha
}^{(n)}\rangle , 
\]
defined by 
\[
\langle \!\langle \Phi ,\varphi \rangle \!\rangle _{\mu }=\langle \!\langle 
\tilde{\Phi},\varphi \rangle \!\rangle _{\tilde{\mu}},\;\varphi \in (%
\mathcal{N})^{1}, 
\]
is in $(\mathcal{N})_{\mu }^{-1}$ and the following relation holds 
\[
\Phi ^{(n)}=\sum_{k+m+l=n}\frac{1}{m!l!}\tilde{\Phi}_{\alpha }^{(k)}\hat{%
\otimes}P_{m}^{\mu ,\alpha }(0)\hat{\otimes}M_{l}^{\tilde{\mu},\alpha }. 
\]
\end{theorem}

\noindent \textbf{Proof. }We can insert formula (\ref{eq2}) in the formula 
\begin{eqnarray*}
&&\sum_{n=0}^{\infty }n!\langle \Phi _{\alpha }^{(n)},\varphi _{\alpha
}^{(n)}\rangle \\
&=&\sum_{k=0}^{\infty }k!\rangle \tilde{\Phi}_{\alpha }^{(k)},\tilde{\varphi}%
_{\alpha }^{(k)}\rangle \\
&=&\sum_{k=0}^{\infty }k!\left\langle \tilde{\Phi}_{\alpha
}^{(k)},\sum_{m,l=0}^{\infty }\frac{(k+m+l)!}{k!m!l!}(P_{m}^{\mu ,\alpha }(0)%
\hat{\otimes}M_{l}^{\tilde{\mu},\alpha },\varphi _{\alpha }^{(k+m+l)})_{%
\mathcal{H}^{\hat{\otimes}(m+l)}}\right\rangle \\
&=&\sum_{k,m,l=0}^{\infty }\frac{(k+m+l)!}{m!l!}\langle \tilde{\Phi}_{\alpha
}^{(k)}\hat{\otimes}P_{m}^{\mu ,\alpha }(0)\hat{\otimes}M_{l}^{\tilde{\mu}%
,\alpha },\varphi _{\alpha }^{(k+m+l)}\rangle \\
&=&\sum_{n=0}^{\infty }n!\left\langle \sum_{k+m+l=n}\frac{1}{m!l!}\tilde{\Phi%
}_{\alpha }^{(k)}\hat{\otimes}P_{m}^{\mu ,\alpha }(0)\hat{\otimes}M_{l}^{%
\tilde{\mu},\alpha },\varphi _{\alpha }^{(n)}\right\rangle ,
\end{eqnarray*}
and compare coefficients again.\hfill $\blacksquare $

\addcontentsline{toc}{section}{References}

\newcommand{\etalchar}[1]{$^{#1}$}

\end{document}